\RequirePackage{ifpdf}
\ifpdf 
\documentclass[pdftex]{sigma}
\else
\documentclass{sigma}
\fi

\numberwithin{equation}{section}

\usepackage[mathscr]{eucal}

\theoremstyle{theorem}
\newtheorem{Theorem}{Theorem}[section]
\newtheorem{Proposition}[Theorem]{Proposition}
\newtheorem{Lemma}[Theorem]{Lemma}
\newtheorem{Corollary}[Theorem]{Corollary}
\theoremstyle{definition}
\newtheorem{Definition}[Theorem]{Definition}
\newtheorem{Example}[Theorem]{Example}
\newtheorem{Remark}[Theorem]{Remark}

\theoremstyle{theorem}
\newtheorem*{claim}{Claim}

\begin{document}

\allowdisplaybreaks

\renewcommand{\PaperNumber}{108}

\FirstPageHeading

\ShortArticleName{T-Systems and Y-Systems for Quantum Af\/f\/inizations}

\ArticleName{T-Systems and Y-Systems for Quantum Af\/f\/inizations\\
of Quantum Kac--Moody Algebras}

\Author{Atsuo KUNIBA~$^\dag$, Tomoki NAKANISHI~$^\ddag$ and Junji SUZUKI~$^\S$}

\AuthorNameForHeading{A. Kuniba, T. Nakanishi and J. Suzuki}

\Address{$^\dag$~Institute of Physics, University of Tokyo, Tokyo, 153-8902, Japan}
\EmailD{\href{mailto:atsuo@gokutan.c.u-tokyo.ac.jp}{atsuo@gokutan.c.u-tokyo.ac.jp}}

\Address{$^\ddag$~Graduate School of Mathematics, Nagoya University, Nagoya, 464-8604, Japan}
\EmailD{\href{mailto:nakanisi@math.nagoya-u.ac.jp}{nakanisi@math.nagoya-u.ac.jp}}

\Address{$^\S$~Department of Physics, Faculty of Science, Shizuoka University, Ohya, 836, Japan}
\EmailD{\href{mailto:sjsuzuk@ipc.shizuoka.ac.jp}{sjsuzuk@ipc.shizuoka.ac.jp}}

\ArticleDates{Received October 05, 2009, in f\/inal form December 16, 2009; Published online December 19, 2009}

\Abstract{The T-systems and Y-systems are
classes of algebraic relations
originally asso\-cia\-ted with quantum af\/f\/ine algebras
and Yangians.
Recently the T-systems were generalized to
quantum af\/f\/inizations
of a wide class of quantum Kac--Moody algebras by
Hernandez.
In this note we introduce
the corresponding Y-systems
and establish a relation
between T and Y-systems.
We also introduce the T and Y-systems
associated with a class of cluster algebras,
which include the
former T and Y-systems of simply laced type
as special cases.}


\Classification{17B37; 13A99}

\medskip

\rightline{\em Dedicated to  Professor Tetsuji Miwa on his 60th birthday}

\section{Introduction}

The T-systems and Y-systems appear in various aspects
for integrable
systems
\cite{Z,KM,KP,KN,RTV,KNS1,KNS2,GT,BLZ,KLWZ,CGT,
DPT,BGOT,DDT,R,GKV,T2}.
Originally, the T-systems are systems of relations
among the Kirillov--Reshetikhin modules \cite{Ki,KR}
in the Grothendieck rings of modules over
quantum af\/f\/ine algebras and Yangians \cite{BR,CP1,KNS1,KS,CP2,FR,N2,Her1,Her4}.
The T and Y-systems are related to each other by certain changes
of variables \cite{KP,KNS1}.

The T and Y-systems are also regarded as
relations among variables for cluster algebras
 \cite{FZ2,FZ4,HL1,DK1,Kel,IIKNS,HL2,N3,DK2}.
This identif\/ication is especially fruitful
in the study of  the periodicity of these systems
\cite{Z,RTV,KNS1,GT,FS,FZ2,FZ3,FZ4,V,Kel,IIKNS}.

The T-systems are generalized by Hernandez \cite{Her3}
to the quantum af\/f\/inizations
of a wide class of quantum Kac--Moody algebras
studied in \cite{D2,VV,Jin,M,N1,Her2}.
In this paper we introduce the corresponding
Y-systems and establish a relation between
T and Y-systems.
We also introduce the
T and Y-systems associated with
a class of  cluster algebras,
which include the former T and Y-systems of simply laced type
as special cases.

It will be interesting to investigate the relation
of the systems discussed here to the birational transformations arising from
the Painlev\'e equations in \cite{NY1,NY2},
and also to the geometric realization of cluster algebras
in \cite{BFZ,GLS}.

The organization of the paper is as follows.
In Section~\ref{sect:qKM} basic def\/initions for quantum Kac--Moody algebras
$U_q(\mathfrak{g})$
and their quantum af\/f\/inizations $U_q(\hat{\mathfrak{g}})$
are recalled.
In Section~\ref{sect:T} the T-systems associated
with the quantum af\/f\/inizations of a class
of quantum Kac--Moody algebras by~\cite{Her3}
are presented.
Based on the result by~\cite{Her3},
the role of the T-system
in the Grothendieck ring
of $U_q(\hat{\mathfrak{g}})$-modules
is given (Corollary \ref{cor:qch1}).
In Section~\ref{sect:Y} we introduce the Y-systems
corresponding to the T-systems in Section~\ref{sect:T},
and establish a relation between them
(Theorem \ref{thm:TtoY1}).
In Section~\ref{sect:rest} we def\/ine the restricted
version of T-systems and Y-systems,
and establish a relation between them
(Theorem \ref{thm:TtoY3}).
In Section~\ref{sect:cluster} we introduce
the T and Y-systems associated with a class of cluster algebras,
which include the restricted T and Y-systems of simply laced
type as special cases.
In particular, the correspondence between
the restricted T and Y-systems of simply laced type
 for the quantum af\/f\/inizations
and cluster algebras is presented
(Corollaries \ref{cor:cluster},
\ref{cor:cluster2}, \ref{cor:cluster3}, and \ref{cor:cluster4}).

\section[Quantum Kac-Moody algebras and their quantum affinizations]{Quantum Kac--Moody algebras\\ and their quantum af\/f\/inizations}
\label{sect:qKM}

In this section, we recall basic def\/initions for
quantum Kac--Moody algebras and their quantum af\/f\/inizations,
following  \cite{Her2,Her3}.
The presentation here is a minimal one.
See \cite{Her2,Her3} for further information and details.

\subsection[Quantum Kac-Moody algebras]{Quantum Kac--Moody algebras}

Let $I=\{1,\dots,r\}$
and let $C=(C_{ij})_{i,j \in I}$ be a
{\em generalized Cartan matrix\/} in
\cite{Ka};
namely, it satisf\/ies
$C_{ij}\in \mathbb{Z}$, $C_{ii}=2$, $C_{ij}\leq 0$ for
any $i\neq j$, and $C_{ij}=0$ if and only if $C_{ji}=0$.
We assume that $C$ is {\em symmetrizable\/},
i.e., there is a diagonal matrix $D=\mathrm{diag}
(d_1,\dots,d_r)$ with $d_i\in \mathbb{N}:=\mathbb{Z}_{> 0}$
such that $B=DC$ is symmetric.
Throughout the paper we assume that there is no common divisor
for $d_1,\dots,d_r$ except for 1.

Let $(\mathfrak{h},\Pi,\Pi^{\vee})$ be a {\em realization\/}
of the Cartan matrix $C$ \cite{Ka};
namely, $\mathfrak{h}$ is a $(2r-\mathrm{rank}\,C)$-dimensional
 $\mathbb{Q}$-vector
space,
and $\Pi=\{\alpha_1,\dots,\alpha_r\}\subset \mathfrak{h}^*$,
$\Pi^{\vee}=\{\alpha_1^{\vee},\dots,\alpha_r^{\vee}\}\subset \mathfrak{h}$
such that $\alpha_j(\alpha_i^{\vee})=C_{ij}$.

Let $q\in \mathbb{C}^{\times}$ be not a root of unity.
We set $q_i = q^{d_i}$ ($i\in I$),
$[k]_q=(q^k-q^{-k})/(q-q^{-1})$,
$[k]_q!=[1]_q[2]_q\cdots [k]_q$,
and ${k \brack r}_q=[k]_q!/[k-r]_q![r]_q!$
($0\leq r\leq k$).

\begin{Definition}[\cite{D1,Jim}]
The {\em quantum Kac--Moody algebra\/} $U_q(\mathfrak{g})$ associated
with $C$ is the $\mathbb{C}$-algebra with generators
$k_h$ ($h\in \mathfrak{h}$), $x_i^{\pm}$ ($i\in I$)
and the following relations:
\begin{gather*}
k_h k_{h'} = k_{h+h'},\qquad
k_0=1,\qquad
k_h x_i^{\pm} k_{-h}
= q^{\pm\alpha_i(h)} x_i^{\pm},\\
x^+_i x^-_j - x^-_j x^+_i = \delta_{ij}
\frac{k_{d_i \alpha_i^{\vee}} - k_{- d_i\alpha_i^\vee}}
{q_i - q_i^{-1}},\\
\sum_{r=0}^{1-C_{ij}}
(-1)^r
{1-C_{ij} \brack{r}}_{q_i}
(x^{\pm}_i)^{1-C_{ij}-r}
x^{\pm}_j
(x^{\pm}_i)^{r}
=0\qquad (i\neq j).
\end{gather*}

\end{Definition}

\subsection{Quantum af\/f\/inizations}

In the following, we use the following formal series
(currents):
\begin{gather*}
x^{\pm}_i(z) =\sum_{r\in \mathbb{Z}} x^{\pm}_{i,r} z^r,\\
\phi^{\pm}_i(z) =\sum_{r\geq 0} \phi^{\pm}_{i,\pm r} z^{\pm r}
:=k_{\pm d_i\alpha_i^{\vee}}
\exp \left(
\pm \big(q-q^{-1}\big)\sum_{r\geq 1} h_{i,\pm r} z^{\pm r}
\right).
\end{gather*}
We also use the {\em formal delta function\/}
$\delta(z)= \sum\limits_{r\in \mathbb{Z}} z^r$.

\begin{Definition}[\cite{D2,Jin,Her3}]
\label{def:Uq}
The {\em quantum affinization\/} (without central elements)
of the quantum Kac--Moody algebra $U_q(\mathfrak{g})$,
denoted by $U_q(\hat{\mathfrak{g}})$,
is the $\mathbb{C}$-algebra with generators
$x^{\pm}_{i,r}$ ($i\in I$, $r\in \mathbb{Z}$),
$k_h$ ($h\in \mathfrak{h}$),
$h_{i,r}$ ($i\in I$, $r\in \mathbb{Z}\setminus\{0\}$)
and the following relations:
\begin{gather}
k_h k_{h'} = k_{h+h'},\qquad
k_0=1,\qquad
k_h \phi^{\pm}_i(z) = \phi^{\pm}_i(z)k_h,\nonumber\\
k_h x^{\pm}_i(z) = q^{\pm \alpha_i(h)}x^{\pm}_i(z)k_h,\nonumber\\
\phi^{+}_i(z)x^{\pm}_j(w)=
\frac{q^{\pm B_{ij}}w-z}{w-q^{\pm B_{ij}}z}
x^{\pm}_j(w)\phi^+_i(z),\nonumber\\
\phi^{-}_i(z)x^{\pm}_j(w)=
\frac{q^{\pm B_{ij}}w-z}{w-q^{\pm B_{ij}}z}
x^{\pm}_j(w)\phi^-_i(z),\nonumber\\
x^+_i(z)x^-_j(w)-x^-_j(w)x^+_i(z)=
\frac{\delta_{ij}}{q_i-q_i^{-1}}
\left(
\delta\left(\frac{w}{z}\right)
\phi^+_i(w)
-
\delta\left(\frac{z}{w}\right)
\phi^-_i(z)
\right),\nonumber\\
(w-q^{\pm B_{ij}}z)x^{\pm}_i(z)x^{\pm}_j(w)=
(q^{\pm B_{ij}}w-z)
x^{\pm}_j(w)x^{\pm}_i(z),\nonumber\\
\sum_{\pi\in \Sigma}
\sum_{k=1}^{1-C_{ij}}
(-1)^k
{1-C_{ij} \brack{k}}_{q_i}
x^{\pm}_i(w_{\pi(1)})\cdots
x^{\pm}_i(w_{\pi(k)})
x^{\pm}_j(z)\nonumber\\
\qquad{}\times
x^{\pm}_i(w_{\pi(k+1)})\cdots
x^{\pm}_i(w_{\pi(1-C_{ij})})=0
\qquad (i\neq j).\label{eq:serre1}
\end{gather}
In \eqref{eq:serre1} $\Sigma$ is the symmetric group
for the set $\{1,\dots,1-C_{ij}\}$.
\end{Definition}

When $C$ is of f\/inite type,
the above $U_q(\hat{\mathfrak{g}})$ is
called  an {\em $($untwisted$)$ quantum affine algebra
$($without central elements$)$} or {\em quantum loop algebra\/};
it is isomorphic to a
subquotient of the quantum Kac--Moody algebra
associated with the (untwisted) af\/f\/ine
extension of $C$
 without derivation
and central elements
 \cite{D2,B}.
(A little confusingly,
 the quantum Kac--Moody algebra associated with $C$ of
  af\/f\/ine type {\em with\/} derivation and central elements
  is also called a quantum af\/f\/ine algebra
and denoted by $U_q(\hat{\mathfrak{g}})$.)

When $C$ is of af\/f\/ine type,
$U_q(\hat{\mathfrak{g}})$
 is called a {\em quantum toroidal algebra\/}
(without central elements).

In general, if $C$ is not of f\/inite type,
$U_q(\hat{\mathfrak{g}})$ is no longer isomorphic to a subquotient of
any quantum Kac--Moody algebra
and has no Hopf algebra structure.


\subsection{The category $\mathrm{Mod}(U_q(\hat{\mathfrak{g}}))$}

Let $U_q(\mathfrak{h})$ be the subalgebra of
$U_q(\hat{\mathfrak{g}})$ generated by $k_h$ ($h\in \mathfrak{h}$).

\begin{Definition}[\cite{Her3}]
Let $\mathrm{Mod}(U_q(\hat{\mathfrak{g}}))$ be the category
of $U_q(\hat{\mathfrak{g}})$-modules satisfying the follo\-wing
properties:
\begin{itemize}
\item[$(a)$] $V$ is $U_q(\mathfrak{h})$-diagonalizable,
i.e., $V=\oplus_{\omega\in \mathfrak{h}^*} V_{\omega}$,
where
\begin{gather*}
V_{\omega}=\{v\in V\mid
k_h v = q^{\omega(h)}v\ \mbox{for any $h\in \mathfrak{h}$}\}.
\end{gather*}
\item[$(b)$] For any
$\omega\in \mathfrak{h}^*$,
$V_{\omega}$ is f\/inite dimensional.
\item[$(c)$] For any $i\in I$ and $\omega\in \mathfrak{h}^*$,
$V_{\omega\pm r\alpha_i}=\{0\}$ for a suf\/f\/iciently large
$r$.
\item[$(d)$] There is a f\/inite number of elements
$\lambda_1,\dots,\lambda_s\in \mathfrak{h}^*$ such that
the weights $\omega\in \mathfrak{h}^*$ of $V$ are in
$\bigcup_{j=1}^s \mathcal{S}(\lambda_j)$,
where $\mathcal{S}(\lambda)=\{ \mu\in \mathfrak{h}^* \mid \mu
\leq \lambda\}$.

\end{itemize}
\end{Definition}
By Condition $(a)$, we restrict our attention to the so called
`type 1' modules.

Let
\begin{gather*}
P=\{ \lambda\in \mathfrak{h}^*
\mid \lambda(\alpha_i^{\vee})\in \mathbb{Z}
\}
\end{gather*}
be the set of integral weights.
\begin{Definition}[\cite{Her2}] \qquad

(1) An {\em $\ell$-weight\/} is a pair $(\lambda,\Psi)$
with $\lambda\in P$ and $\Psi=(\Psi^{\pm}_{i,\pm r})_{i\in I,r\geq 0}$
such that $\Psi^{\pm}_{i,\pm r}\in \mathbb{C}$ and
$\Psi^{\pm}_{i,0}=q_i^{\pm\lambda(\alpha^{\vee}_i)}$.

(2) A $U_q(\hat{\mathfrak{g}})$-module $V$ is
{\em of $\ell$-highest weight\/}
if there is some $v\in V$
and $\ell$-weight $(\lambda,\Psi)$
 such that
$x^+_{i,r}v=0$, $k_h v=q^{\lambda(h)} v$,
$\phi^{\pm}_{i,\pm r}v = \Psi^{\pm}_{i,\pm r} v$,
and $U_q(\hat{\mathfrak{g}})v=V$.
Such $v$ and
$(\lambda,\Psi)$ are called
a~{\em highest weight vector\/} and
the  {\em $\ell$-highest weight\/} of $V$, respectively.
\end{Definition}

By the standard argument using Verma modules, one can show that
for any $\ell$-weight $(\lambda,\Psi)$, there is a unique
simple $\ell$-highest weight module $L(\lambda,\Psi)$ with
$\ell$-highest weight $(\lambda,\Psi)$ \cite{Her2}.

The following theorem is a generalization of the well-known
classif\/ication of the simple f\/inite-dimensional
modules of the quantum af\/f\/ine algebras by \cite{CP1,CP2}.

\begin{Theorem}[\cite{M,N1,Her2}]
We have $L(\lambda,\Psi)\in
 \mathrm{Mod}(U_q(\hat{\mathfrak{g}}))$ if and only if
there is an $I$-tuple of polynomials $(P_i(u))_{i\in I}$,
$P_i(u)\in \mathbb{C}[u]$ with $P_i(0)=1$ such that
\begin{gather*}
\sum_{m\geq 0} \Psi^{\pm}_{i,\pm m} z^{\pm m}
=q_i^{\mathrm{deg}\,P_i}
\frac{P_i(zq_i^{-1})}{P_i(zq_i)},
\end{gather*}
where the equality is in $\mathbb{C}[[z^{\pm1}]]$.
\end{Theorem}

We call $(P_i(u))_{i\in I}$
the {\em Drinfeld polynomials\/} of $L(\lambda,\Psi)$.
In the case of quantum af\/f\/ine algebras,
$\lambda$ is also completely determined by the Drinfeld polynomials
by the condition $\lambda(\alpha_i^{\vee})=\mathrm{deg}\,P_i$.
This is not so in general.

Let  $\Lambda_i\in \mathfrak{h}^*$ ($i\in I$)
be the {\em fundamental weights
of $U_q(\mathfrak{g})$} satisfying
 $\Lambda_i(\alpha_j^{\vee})=\delta_{ij}$.

\begin{Example}
\label{example:fund}
 For the  following choices of $(\lambda,\Psi)$,
$L(\lambda,\Psi)\in
 \mathrm{Mod}(U_q(\hat{\mathfrak{g}}))$ is
called a {\em fundamental module} \cite{Her3}.

$(a)$ For any $i\in I$ and $\alpha\in \mathbb{C}^{\times}$,
set  $\lambda=\Lambda_i$, $P_i(u)=1-\alpha u$,
and $P_j(u)=1$ ($j\in I, j\neq i$).

$(b)$ Choose any $\lambda$ satisfying $(\lambda,\alpha_i^{\vee})=0$
 ($i\in I$) and  also set  $P_i(u)=1$ ($i\in I$).
The corresponding module $L(\lambda,\Psi)$
  is written  as $L(\lambda)$.
The module $L(\lambda)$ is one-dimensional;
it is trivial in the case of the quantum af\/f\/ine algebras.
\end{Example}

\subsection[Kirillov-Reshetikhin modules]{Kirillov--Reshetikhin modules}

The following is a generalization of the Kirillov--Reshetikhin modules
of the quantum af\/f\/ine algebras studied by
\cite{Ki,KR,BR,CP2,KNS1,CP3}.

\begin{Definition}[\cite{Her3}]
\label{def:KR1}
 For any $i\in I$, $m\in\mathbb{N}$,
and $\alpha\in \mathbb{C}^{\times}$,
set the polynomials
$(P_j(u))_{j\in I}$ as
\begin{gather}
\label{eq:KR1}
P_i(u)=\big(1-\alpha q_i^{m-1}u\big)\big(1-\alpha q_i^{m-3}u\big)
\cdots \big(1-\alpha q_i^{1-m}u\big)
\end{gather}
and $P_j(u)=1$ for any $j\neq i$.
The corresponding module $L(m\Lambda_i,\Psi)$ is called
a {\em Kirillov--Reshetikhin module\/}
and denoted by $W^{(i)}_{m,\alpha}$.
\end{Definition}

\begin{Remark}
We slightly shift the def\/inition of the polynomial \eqref{eq:KR1}
for $W^{(i)}_{m,\alpha}$ in \cite{Her3} in order to make
the identif\/ication to the forthcoming T-systems a little simpler.
\end{Remark}

\section{T-systems}
\label{sect:T}

\subsection{T-systems}
Throughout  Sections \ref{sect:T}--\ref{sect:rest},
 we restrict our attention
to a symmetrizable generalized Cartan matrix~$C$
satisfying the following condition due to Hernandez
\cite{Her3}:
\begin{gather}
\label{eq:Ccond1}
\mbox{If $C_{ij}< -1$, then $d_i=-C_{ji}=1$,}
\end{gather}
where $D=\mathrm{diag}(d_1,\dots,d_r)$
is the diagonal matrix symmetrizing $C$.
In this paper, we say that a
generalized Cartan matrix $C$ is {\em tamely laced\/} if
it is symmetrizable and satisf\/ies the condition~\eqref{eq:Ccond1}.

As usual,
 we say that a generalized Cartan matrix $C$ is {\em simply laced\/} if
$C_{ij}=0$ or $-1$ for any $i\neq j$.
If $C$ is simply laced, then it is symmetric,
$d_a=1$ for any $a\in I$,
and it is tamely laced.

With a tamely laced generalized Cartan matrix $C$,
we  associate a {\em Dynkin diagram\/} in the standard way \cite{Ka}:
For any pair $i\neq j \in I$ with $C_{ij}<0$,
the vertices $i$ and $j$ are connected by
 $\max\{|C_{ij}|,|C_{ji}|\}$ lines,
and the lines are equipped with an arrow from $j$ to $i$
 if $C_{ij}<-1$.
Note that the condition \eqref{eq:Ccond1} means

\begin{itemize}\itemsep=0pt
\item[$(i)$] the vertices $i$ and $j$ are not connected
if $d_i,d_j>1$ and $d_i\neq d_j$,

\item[$(ii)$]
the vertices $i$ and $j$ are connected by $d_i$ lines with
 an arrow from $i$ to $j$
or not connected if $d_i>1$ and $d_j=1$,

\item[$(iii)$]
the vertices $i$ and $j$ are connected by a single line
or not connected if $d_i=d_j$.
\end{itemize}

\begin{Example}\qquad
$(1)$ Any Cartan matrix of f\/inite or af\/f\/ine type is
tamely laced except for types~$A^{(1)}_1$ and $A^{(2)}_{2\ell}$.

$(2)$ The following generalized Cartan matrix $C$
is tamely laced:
\begin{gather*}
C=
\begin{pmatrix}
2 & -1 & 0 & 0\\
-3& 2 & -2 & -2\\
0 & -1 & 2 & -1\\
0 & -1 & -1 & 2\\
\end{pmatrix},
\qquad
D=
\begin{pmatrix}
3 & 0 & 0 & 0\\
0& 1 & 0 & 0\\
0 & 0 & 2 & 0\\
0 & 0 & 0 & 2\\
\end{pmatrix}.
\end{gather*}
The corresponding Dynkin diagram is
\begin{gather*}
\includegraphics{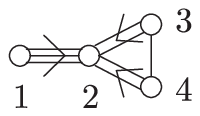}
\end{gather*}
\end{Example}

For a tamely laced generalized Cartan matrix $C$, we set  an integer $t$ by
\begin{gather*}
t=\mathrm{lcm}(d_1,\dots,d_r).
\end{gather*}
For $a,b\in I$, we write $a\sim b$ if
$C_{ab}<0$, i.e., $a$ and $b$ are adjacent
in the corresponding Dynkin diagram.

Let $U$ be either $\frac{1}{t}\mathbb{Z}$,
 the complex plane $\mathbb{C}$, or
the cylinder $\mathbb{C}_{\xi}:= \mathbb{C}/
(2\pi \sqrt{-1}/\xi) \mathbb{Z}$ for some
 $\xi\in \mathbb{C}
\setminus 2\pi \sqrt{-1}\mathbb{Q}$,
depending on the situation under consideration.
The following is a generalization of the T-systems associated
with the quantum af\/f\/ine algebras~\cite{KNS1}.

\begin{Definition}[\cite{Her3}]
For a tamely laced generalized Cartan matrix $C$,
the {\em unrestricted T-system $\mathbb{T}(C)$ associated with $C$}
is the following system of relations for
a family of variables
 $T=\{T^{(a)}_m(u) \mid a\in I, m\in \mathbb{N},
u\in U \}$,
\begin{gather}
\label{eq:T1}
T^{(a)}_m\left(u-\frac{d_a}{t}\right)
T^{(a)}_m\left(u+\frac{d_a}{t}\right)
 =
T^{(a)}_{m-1}(u)T^{(a)}_{m+1}(u)
+
\prod_{b: b\sim a}
T^{(b)}_{\frac{d_a}{d_b}m}(u)
\qquad\mbox{if \ \ $d_a>1$},
\\
\label{eq:T2}
T^{(a)}_m\left(u-\frac{d_a}{t}\right)
T^{(a)}_m\left(u+\frac{d_a}{t}\right)
 =
T^{(a)}_{m-1}(u)T^{(a)}_{m+1}(u)
+
\prod_{b: b\sim a}
S^{(b)}_{m}(u)
\qquad\mbox{if \ \ $d_a=1$},
\end{gather}
where
$T^{(a)}_0 (u)= 1$ if they occur
in the right hand sides in the relations.
The symbol $S^{(b)}_{m}(u)$ is def\/ined by
\begin{gather}
\label{eq:Sbm}
S^{(b)}_{m}(u)=
{\prod_{k=1}^{d_b}}T^{(b)}_{1+E\left[\frac{m-k}{d_b}\right]}
\left(u+\frac{1}{t}\left(2k-1-m+E\left[\frac{m-k}{d_b}\right]d_b\right)
\right),
\end{gather}
and $E[x]$ ($x\in \mathbb{Q}$) denotes the largest integer
not exceeding $x$.
\end{Definition}

\begin{Remark}\qquad

1. This is a slightly reduced version
of the T-systems in \cite[Theorem 6.10]{Her3}.
See Remark \ref{rem:T1}.
The same system was also studied by \cite{T1}
when $C$ is of af\/f\/ine type
in view of a~generalization of discrete Toda f\/ield equations.

2. More explicitly, $S^{(b)}_m(u)$ is written as follows:
For $0\leq j < d_b$,
\begin{gather*}
S^{(b)}_{d_bm+j}(u)
=
\left\{
\prod_{k=1}^j
 T^{(b)}_{m+1}
\left(u+\frac{1}{t}(j+1-2k)\right)
\right\}\!
\left\{
\prod_{k=1}^{d_b-j}
 T^{(b)}_{m}\left(u+\frac{1}{t}(d_b-j+1-2k)\right)
\right\}\!.\!\!\!
\end{gather*}
For example, for $d_b=1$,
\begin{gather*}
S^{(b)}_{m}(u)=T^{(b)}_m(u),
\end{gather*}
for $d_b=2$,
\begin{gather*}
S^{(b)}_{2m}(u)
 =
T^{(b)}_{m}\left(u- \frac{1}{t}\right)
T^{(b)}_{m}\left(u+\frac{1}{t}\right),\\
S^{(b)}_{2m+1}(u)
 =
T^{(b)}_{m+1}(u)
T^{(b)}_{m}(u),
\end{gather*}
for $d_b=3$,
\begin{gather*}
S^{(b)}_{3m}(u)
 =
T^{(b)}_{m}\left(u- \frac{2}{t}\right)
T^{(b)}_{m}(u)
T^{(b)}_{m}\left(u+\frac{2}{t}\right),\\
S^{(b)}_{3m+1}(u)
 =
T^{(b)}_{m+1}(u)
T^{(b)}_{m}\left(u- \frac{1}{t}\right)
T^{(b)}_{m}\left(u+\frac{1}{t}\right),\\
S^{(b)}_{3m+2}(u)
 =
T^{(b)}_{m+1}\left(u- \frac{1}{t}\right)
T^{(b)}_{m+1}\left(u+ \frac{1}{t}\right)
T^{(b)}_{m}(u),
\end{gather*}
and so on.

3. The second terms in the right hand sides of
\eqref{eq:T1} and \eqref{eq:T2} can be written in a unif\/ied way
as follows \cite{Her3}:
\begin{gather*}
\prod_{b:b\sim a}
\prod_{k=1}^{-C_{ab}}
T^{(b)}_{-C_{ba}+E\left[\frac{d_a(m-k)}{d_b}\right]}
\left(u+\frac{d_b}{t}\left(\frac{-2k+1}{C_{ab}}-C_{ba}
+E\left[\frac{d_a(m-k)}{d_b}\right]-1\right)
-\frac{d_am}{t}
\right).
\end{gather*}
\end{Remark}

\begin{Definition}
\label{defn:TC}
Let $\EuScript{T}(C)$
be the commutative ring over $\mathbb{Z}$  with generators
$T^{(a)}_m(u)^{\pm 1}$ ($a\in I$, $m\in \mathbb{N}$,
$u\in U $)
and the relations $\mathbb{T}(C)$.
(Here we also assume the relation
$T^{(a)}_m(u)T^{(a)}_m(u)^{-1}=1$ implicitly.
We do not repeat this remark in the forthcoming
similar def\/initions.)
Also, let $\EuScript{T}^{\circ}(C)$
be the subring of $\EuScript{T}(C)$
generated by
$T^{(a)}_m(u)$ ($a\in I$, $m\in \mathbb{N}$,
$u\in U $).
\end{Definition}

\subsection{T-system and Grothendieck ring}
\label{subsect:rep}

Let $C$ continue to be a tamely laced generalized Cartan matrix.
The T-system $\mathbb{T}(C)$ is a family of
relations in the Grothendieck ring of modules
of $U_q(\hat{\mathfrak{g}})$ as explained below.

We recall  facts on $\mathrm{Mod}(U_q(\hat{\mathfrak{g}}))$
in \cite{Her3}.
\begin{itemize}\itemsep=0pt
\item[1.] For a pair of $\ell$-highest weight modules $V_1, V_2\in
\mathrm{Mod}(U_q(\hat{\mathfrak{g}}))$,
there is an $\ell$-highest weight module $V_1*_f V_2\in
\mathrm{Mod}(U_q(\hat{\mathfrak{g}}))$
called the {\em fusion product\/}.
It is def\/ined by using the $u$-deformation of the
Drinfeld coproduct and the  specialization at $u=1$.

\item[2.] Any $\ell$-highest weight module in
$\mathrm{Mod}(U_q(\hat{\mathfrak{g}}))$ has a f\/inite composition.
(Here, the condition~\eqref{eq:Ccond1} for $C$ is
used essentially in \cite{Her3}.)

\item[3.] If $V_1$, $V_2\in
\mathrm{Mod}(U_q(\hat{\mathfrak{g}}))$ have f\/inite compositions,
then $V_1 *_f V_2$ also has a f\/inite composi\-tion.
\end{itemize}
Therefore, the {\em Grothendieck ring
$R(C)$ of the modules   in
$\mathrm{Mod}(U_q(\hat{\mathfrak{g}}))$
having finite compositions\/}
is well def\/ined, where the product is given by~$*_f$.

Let $R'(C)$ be the quotient
ring  of $R(C)$ by the ideal generated
by all $L(\lambda,\Psi)-L(\lambda',\Psi)$'s.
In other words, we regard modules in $R(C)$
as modules of the subalgebra
of $U_q(\hat{\mathfrak{g}})$ generated by~$x^{\pm}_{i,r}$ ($i\in I$, $r\in \mathbb{Z}$),
$k_{\pm d_i \alpha^{\vee}_i}$ ($i\in I$),
$h_{i,r}$ ($i\in I$, $r\in \mathbb{Z}\setminus\{0\}$)
in Def\/inition \ref{def:Uq}.

\begin{Proposition}
\label{prop:Rgen}
The ring $R'(C)$ is freely generated by
 the fundamental modules $L(\Lambda_i,\Psi)$
in Example $\ref{example:fund}~(a)$.
\end{Proposition}

\begin{proof}
It follows from \cite[Corollary 4.9]{Her3} that
$R'(C)$ is generated by
 the fundamental modules $L(\Lambda_i,\Psi)$.
The $q$-character morphism $\chi_q$ def\/ined in~\cite{Her3}
induces an injective ring homomorphism
$\chi_q':R'(C)\rightarrow \mathbb{Z}[Y_{i,\alpha}^{\pm1}]_{i\in I, \alpha\in
\mathbb{C}^{\times}}$.
Furthermore,
for $L(\Lambda_i,\Psi)$ with $P_i(u)=1-\alpha u$,
the highest term of $\chi_q'(L(\Lambda_i,\Psi))$
is $Y_{i,\alpha}$.
Therefore, the algebraic independence of
$L(\Lambda_i,\Psi)$'s follows from that of~$Y_{i,\alpha}$'s.
\end{proof}

We set $\mathbb{C}_{t\log q}:=
\mathbb{C}/(2\pi \sqrt{-1}/(t\log q))\mathbb{Z}$,
and introduce alternative notation $W^{(a)}_m(u)$
 ($a\in I$, $m\in \mathbb{N}$, $u\in \mathbb{C}_{t\log q}$)
for the Kirillov--Reshetikhin module $W^{(a)}_{m,q^{tu}}$
in Def\/inition~\ref{def:KR1}.

In terms of the Kirillov--Reshetikhin
modules, the structure of $R'(C)$ is described
as follows:
\begin{Theorem}
\label{thm:qch1}
Let
$W=\{W^{(a)}_m(u)
\mid
a\in I, m\in \mathbb{N},
u\in \mathbb{C}_{t\log q}
\}$
 be the family
of the Kirillov--Reshetikhin modules in $R'(C)$.
Let $T$ and $\mathbb{T}(C)$ be the ones in
Definition $\ref{defn:TC}$
with $U=\mathbb{C}_{t\log q}$.
Then,

$(1)$ The family $W$ generates the ring $R'(C)$.

$(2)$ {\rm (\cite{Her3})}
The family $W$ satisfies the T-system $\mathbb{T}(C)$
in $R'(C)$
by replacing $T^{(a)}_m(u)$  in $\mathbb{T}(C)$
 with $W^{(a)}_m(u)$.

$(3)$ For any polynomial $P(T)\in \mathbb{Z}[T]$,
the relation $P(W)=0$ holds
 in $R'(C)$
if and only if there is a nonzero monomial $M(T)\in \mathbb{Z}[T]$
such that  $M(T)P(T) \in I(\mathbb{T}(C))$,
where $I(\mathbb{T}(C))$ is the ideal of $\mathbb{Z}[T]$
generated by the relations in $\mathbb{T}(C)$.
\end{Theorem}

\begin{proof}
(1) By Proposition \ref{prop:Rgen},
$R'(C)$ is generated by
$L(\Lambda_i,\Psi)$'s, which belong to $W$.

(2) This is due to  \cite[Theorem 6.10]{Her3}.

(3) The proof is the same with that of \cite[Theorem 2.8]{IIKNS}
by generalizing the {\em height\/} of $T^{(a)}_m(u)$ therein
as $\mathrm{ht}\, T^{(a)}_m(u):= d_a (m-1)+1$.
\end{proof}

\begin{Remark}
\label{rem:T1}
In \cite{Her3} the T-system is considered in $R(C)$
including fundamental modules $L(\lambda)$ of Example
\ref{example:fund}~$(b)$.
\end{Remark}

As a corollary, we have a generalization of \cite[Corollary 2.9]{IIKNS}
for the quantum af\/f\/ine algebras:
\begin{Corollary}
\label{cor:qch1}
The ring
 $\EuScript{T}^{\circ}(C)$
with $U=\mathbb{C}_{t\log q}$
is isomorphic to $R'(C)$
by the correspondence $T^{(a)}_m(u)\mapsto W^{(a)}_m(u)$.
\end{Corollary}

\section{Y-systems}
\label{sect:Y}

\subsection{Y-systems}
\begin{Definition}
For a tamely laced generalized Cartan matrix $C$,
the {\em unrestricted Y-system $\mathbb{Y}(C)$ associated with~$C$}
is the following system of relations for
a family of variables
 $Y=\{Y^{(a)}_m(u) \mid a\in I, m\in \mathbb{N},
u\in U \}$,
where
$Y^{(a)}_0 (u)^{-1}= 0$ if they occur
in the right hand sides in the relations:
\begin{gather}
\label{eq:Y1}
Y^{(a)}_m\left(u-\frac{d_a}{t}\right)
Y^{(a)}_m\left(u+\frac{d_a}{t}\right)
=
\frac{
{ \prod\limits_{b:b\sim a}}
Z^{(b)}_{\frac{d_a}{d_b},m}(u)
}
{
(1+Y^{(a)}_{m-1}(u)^{-1})(1+Y^{(a)}_{m+1}(u)^{-1})}
\qquad\mbox{if \ \ $d_a>1$},
\\
\label{eq:Y2}
Y^{(a)}_m\left(u-\frac{d_a}{t}\right)
Y^{(a)}_m\left(u+\frac{d_a}{t}\right)
 =
\frac{
{\prod\limits_{b:b\sim a}}
\big(1+Y^{(b)}_{\frac{m}{d_b}}(u)\big)
}
{
(1+Y^{(a)}_{m-1}(u)^{-1})(1+Y^{(a)}_{m+1}(u)^{-1})}
\qquad\mbox{if \ \ $d_a=1$},
\end{gather}
where for $p\in \mathbb{N}$
\begin{gather*}
Z^{(b)}_{p,m}(u)=
{\prod_{j=-p+1}^{p-1}}
\left\{
{\prod_{k=1}^{p-|j|}}
\left(
1+Y^{(b)}_{pm+j}\left(
u+\frac{1}{t}(p-|j|+1-2k)
\right)
\right)
\right\},
\end{gather*}
and $Y^{(b)}_{\frac{m}{d_b}}(u)=0$ in \eqref{eq:Y2}
if $\frac{m}{d_b}\not\in \mathbb{N}$.
\end{Definition}

\begin{Remark}\qquad

1. The Y-systems here are formally in the same form
as the ones for the quantum af\/f\/ine algebras \cite{KN}.
However, $p$ for $Z^{(b)}_{p,m}(u)$
here may be greater than 3.

2. In the right hand side of \eqref{eq:Y1},
$\frac{d_a}{d_b}$ is either 1 or $d_a$
due to \eqref{eq:Ccond1}.
The term
 $Z^{(b)}_{p,m}(u)$ is written more explicitly as follows:
for $p=1$,
\begin{gather*}
Z^{(b)}_{1,m}(u)
=1+Y^{(b)}_{m}(u),
\end{gather*}
for $p=2$,
\begin{gather*}
Z^{(b)}_{2,m}(u)
=\big(1+Y^{(b)}_{2m-1}(u)\big)
\left(1+Y^{(b)}_{2m}\left(u-\frac{1}{t}\right)\right)\left(1+Y^{(b)}_{2m}\left(u+\frac{1}{t}\right)\right)
\big(1+Y^{(b)}_{2m+1}(u)\big),
\end{gather*}
for $p=3$,
\begin{gather*}
Z^{(b)}_{3,m}(u)
=
\big(1+Y^{(b)}_{3m-2}(u)\big)
\left(1+Y^{(b)}_{3m-1}\left(u-\frac{1}{t}\right)\right)\left(1+Y^{(b)}_{3m-1}\left(u+\frac{1}{t}\right)\right)\\
\phantom{Z^{(b)}_{3,m}(u)=}{}
\times \left(1+Y^{(b)}_{3m}\left(u-\frac{2}{t}\right)\right)
\big(1+Y^{(b)}_{3m}(u)\big)
\left(1+Y^{(b)}_{3m}\left(u+\frac{2}{t}\right)\right)\\
\phantom{Z^{(b)}_{3,m}(u)=}{}
\times \left(1+Y^{(b)}_{3m+1}\left(u-\frac{1}{t}\right)\right)\left(1+Y^{(b)}_{3m+1}\left(u+\frac{1}{t}\right)\right)
\big(1+Y^{(b)}_{3m+2}(u)\big),
\end{gather*}
and so on.
There are $p^2$ factors in $Z^{(b)}_{p,m}(u)$.
\end{Remark}

\subsection{Relation between T and Y-systems}

Let us write both the relations \eqref{eq:T1}
and \eqref{eq:T2} in $\mathbb{T}(C)$ in a
 unif\/ied manner
\begin{gather}
T^{(a)}_m\left(u-\frac{d_a}{t}\right)T^{(a)}_m\left(u+\frac{d_a}{t}\right)=
T^{(a)}_{m-1}(u)T^{(a)}_{m+1}(u)
+M^{(a)}_m(u)\nonumber\\
\hphantom{T^{(a)}_m\left(u-\frac{d_a}{t}\right)T^{(a)}_m\left(u+\frac{d_a}{t}\right)}{} =
T^{(a)}_{m-1}(u)T^{(a)}_{m+1}(u)
+\prod_{(b,k,v)}T^{(b)}_k(v)^{G(b,k,v;\,a,m,u)},\label{eq:ta1}
\end{gather}
where $M^{(a)}_m(u)$ is the second term of the right
hand side of each relation.
Def\/ine the transposition $^t G(b,k,v;\,a,m,u)=G(a,m,u;\,b,k,v)$.

\begin{Theorem} The Y-system $\mathbb{Y}(C)$ is
written as
\begin{gather*}
Y^{(a)}_m\left(u-\frac{d_a}{t}\right)Y^{(a)}_m\left(u+\frac{d_a}{t}\right)
=
\frac{
\prod_{(b,k,v)}\big(1+Y^{(b)}_k(v)\big)^{^t\!G(b,k,v;\,a,m,u)}
}
{
\big(1+Y^{(a)}_{m-1}(u)^{-1}\big)\big(1+Y^{(a)}_{m+1}(u)^{-1}\big)
}.
\end{gather*}
\end{Theorem}
\begin{proof}
This can be proved by case check for $d_a>1$ and $d_a=1$.
\end{proof}

For any commutative ring $R$ over $\mathbb{Z}$ with identity element,
let $R^{\times}$ denote the group of all the invertible elements of $R$.

\begin{Theorem}
\label{thm:TtoY1}
Let $R$ be any commutative ring over $\mathbb{Z}$ with identity element.

$(1)$ For any family
$T=\{T^{(a)}_m(u)\in R^{\times} \mid a\in I, m\in \mathbb{N}, u\in U\}$
satisfying $\mathbb{T}(C)$ in $R$,
define a~family
$Y=\{Y^{(a)}_m(u)\in R^{\times} \mid a\in I, m\in \mathbb{N}, u\in U\}$ by
\begin{gather}
\label{eq:TtoY1}
Y^{(a)}_m(u)=\frac{M^{(a)}_m(u)}
{T^{(a)}_{m-1}(u)T^{(a)}_{m+1}(u)},
\end{gather}
where $T^{(a)}_0(u)=1$.
Then,
\begin{gather}
\label{eq:TtoY2}
1+Y^{(a)}_m(u) =
\frac{T^{(a)}_{m}\bigl(u- \frac{d_a}{t}\bigr)
T^{(a)}_{m}\left(u+ \frac{d_a}{t}\right)}
{T^{(a)}_{m-1}(u)T^{(a)}_{m+1}(u)},\\
\label{eq:TtoY3}
1+Y^{(a)}_m(u)^{-1} =
\frac{T^{(a)}_{m}\bigl(u- \frac{d_a}{t}\bigr)
T^{(a)}_{m}\left(u+ \frac{d_a}{t}\right)}
{M^{(a)}_m(u)}.
\end{gather}
Furthermore,
$Y$ satisfies $\mathbb{Y}(C)$ in $R$.

$(2)$ Conversely,
for any family
$Y=\{Y^{(a)}_m(u)\in R^{\times} \mid a\in I, m\in \mathbb{N}, u\in U\}$
 satisfying $\mathbb{Y}(C)$
with $1+Y^{(a)}_m(u)^{\pm1}\in R^{\times}$,
there is a $($not unique$)$ family
$T=\{T^{(a)}_m(u)\in R^{\times} \mid a\in I, m\in \mathbb{N},$ $u\in U\}$
 satisfying $\mathbb{T}(C)$ such that
$Y^{(a)}_m(u)$ is given by~\eqref{eq:TtoY1}.
\end{Theorem}

\begin{proof}
(1)  Equations~\eqref{eq:TtoY2} and
\eqref{eq:TtoY3} follow from  \eqref{eq:ta1} and \eqref{eq:TtoY1}.
We show that $Y$ satis\-f\/ies~$\mathbb{Y}(C)$.
For $d_a>1$,
by \eqref{eq:TtoY1}--\eqref{eq:TtoY3},
the relation \eqref{eq:Y1} reduces to the
following identity:
For any $p\in \mathbb{N}$,
\begin{gather}
\label{eq:Tidentity1}
\frac{
T^{(b)}_{pm}\left(u-\frac{p}{t}\right)
T^{(b)}_{pm}\left(u+\frac{p}{t}\right)
}
{
T^{(b)}_{p(m-1)}(u)
T^{(b)}_{p(m+1)}(u)
}
 =
{ \prod_{j=-p+1}^{p-1}}
\left\{
{ \prod_{k=1}^{p-|j|}}
\frac{
T^{(b)}_{pm+j}\left(\tilde{u}-\frac{1}{t}\right)
T^{(b)}_{pm+j}\left(\tilde{u}+\frac{1}{t}\right)
}
{
T^{(b)}_{pm+j-1}(\tilde{u})
T^{(b)}_{pm+j+1}(\tilde{u})
}
\right\},
\end{gather}
where $\tilde{u}=u+\frac{1}{t}(p-|j|+1-2k)$.
This is easily proved without
using $\mathbb{T}(C)$.
Similarly, for $d_a=1$,
the relation \eqref{eq:Y2} reduces to the
following identity:
\begin{gather}
\label{eq:Tidentity2}
\frac{
S^{(b)}_{m}\left(u-\frac{1}{t}\right)
S^{(b)}_{m}\left(u+\frac{1}{t}\right)
}
{
S^{(b)}_{m-1}(u)
S^{(b)}_{m+1}(u)
}
 =
\begin{cases}
\displaystyle
\frac{
T^{(b)}_{\frac{m}{d_b}}\bigl(u-\frac{d_b}{t}\bigr)
T^{(b)}_{\frac{m}{d_b}}\bigl(u+\frac{d_b}{t}\bigr)
}
{
T^{(b)}_{\frac{m}{d_b}-1}(u)
T^{(b)}_{\frac{m}{d_b}+1}(u)
},
& \displaystyle \frac{m}{d_b}\in \mathbb{N},\\
1,& \mbox{otherwise},
\end{cases}
\end{gather}
where $S^{(b)}_m(u)$ is def\/ined in
\eqref{eq:Sbm}.
Again, this is easily proved without
using $\mathbb{T}(C)$.

(2)
We modify the proof in the case of  quantum af\/f\/ine algebras
  \cite[Theorem 2.12]{IIKNS}
so that it is  applicable to the present situation.
Here, we concentrate on the case $U=\frac{1}{t}\mathbb{Z}$.
The modif\/ication of the proof
 for the
other cases $U=\mathbb{C}$ and $\mathbb{C}/(2\pi
\sqrt{-1}/\xi)\mathbb{Z}$ is straightforward.

{\bf Case 1.} {\em When $C$ is simply laced.}
Suppose that $C$ is simply laced.
Thus, $d_a=1$ for any $a\in I$ and $t=1$.
For any $Y$ satisfying $\mathbb{Y}(C)$,
we construct a desired family $T$ in the following
three steps:

{\em Step 1.} Choose arbitrarily
 $T^{(a)}_1(-1), T^{(a)}_1(0), \in R^{\times}$ ($a\in I$).

{\em Step 2.}
Def\/ine $T^{(a)}_1(-2)$, $T^{(a)}_1(1)$ ($a\in I$)
by
\begin{gather}
\label{eq:ty4}
T^{(a)}_{1}(u\pm 1) =
\big(1+Y^{(a)}_1(u)^{-1}\big)
\frac{M^{(a)}_1(u)}
{T^{(a)}_{1}(u\mp 1)}.
\end{gather}
Repeat it
and def\/ine $T^{(a)}_1(u)$ ($a\in I$) for the rest of  $u\in \mathbb{Z}$
by \eqref{eq:ty4}.

{\em Step 3.}
Def\/ine $T^{(a)}_m(u)$ ($a \in I, u\in \mathbb{Z}$) for
$m\geq 2$ by
\begin{gather}
\label{eq:ty5}
T^{(a)}_{m+1}(u) =
\frac{1}{1+Y^{(a)}_m(u)}
\frac{T^{(a)}_m(u-1)T^{(a)}_m(u+1)}
{T^{(a)}_{m-1}(u)},
\end{gather}
where $T^{(a)}_0(u)=1$.

\begin{claim}
The family $T$ defined above satisfies the
following relations in $R$ for any $a\in I$, $m\in \mathbb{N}$, $u\in
\mathbb{Z}$:
\begin{gather}
  \label{eq:T-Y1}
  Y_m^{(a)}(u)
   =
  \frac{M_m^{(a)}(u)}
  {T_{m-1}^{(a)}(u) T_{m+1}^{(a)}(u)},\\
  \label{eq:T-Y2}
  1+Y_m^{(a)}(u)
   =
  \frac{T_m^{(a)}(u-1)
 T_m^{(a)}(u+1)}
  {T_{m-1}^{(a)}(u) T_{m+1}^{(a)}(u)},\\
  \label{eq:T-Y3}
  1+Y_m^{(a)}(u)^{-1}
   =
  \frac{T_m^{(a)}(u-1)
 T_m^{(a)}(u+1)}
       {M_m^{(a)}(u)}.
\end{gather}
\end{claim}

\noindent
{\em Proof of Claim.} \eqref{eq:T-Y2} holds by \eqref{eq:ty5}.
\eqref{eq:T-Y1} and
\eqref{eq:T-Y3} are equivalent under \eqref{eq:T-Y2};
furthermore, \eqref{eq:T-Y3} holds for
$m=1$ by \eqref{eq:ty4}.
So it is enough to show \eqref{eq:T-Y1} for $m\geq 2$
by induction on~$m$. In fact,
\begin{gather*}
  \frac{M_m^{(a)}(u)}
  {T_{m-1}^{(a)}(u)
 T_{m+1}^{(a)}(u)}
 \overset{\text{by \eqref{eq:T-Y2}}}{=}
  \frac  {(1+Y_{m-2}^{(a)}(u)) T_{m-3}^{(a)}(u)}
{T_{m-2}^{(a)}(u-1)T_{m-2}^{(a)}(u+1)}
  \frac  {(1+Y_{m}^{(a)}(u)) T_{m-1}^{(a)}(u)}
{T_{m}^{(a)}(u-1)T_{m}^{(a)}(u+1)}
{M^{(a)}_{m}(u)}
\\
\overset{\mbox{\scriptsize $\begin{array}{@{}c@{}}{\rm by \ induction}\\ {\rm hypothesis}\end{array}$}}{=}\!\!
  \big(1+Y_{m-2}^{(a)}(u)\big) \big(1+Y_{m}^{(a)}(u)\big)
  \frac{Y_{m-1}^{(a)}(u-1)Y_{m-1}^{(a)}(u+1)}
  { Y_{m-2}^{(a)}(u)}
\frac
{
 M^{(a)}_{m-2}(u)
 M^{(a)}_{m}(u)
}
{ M^{(a)}_{m-1}(u-1)
 M^{(a)}_{m-1}(u+1)}\!
\\
\overset{\text{by \eqref{eq:T-Y2} and $\mathbb{Y}(C)$}}{=}
   Y_{m}^{(a)}(u).
\end{gather*}
This ends the proof of Claim.

Now,
 taking the inverse sum of
(\ref{eq:T-Y2}) and (\ref{eq:T-Y3}),
we obtain (\ref{eq:ta1}).
Therefore, $T$ satisf\/ies the desired properties.

{\bf Case 2.} {\em When $C$ is nonsimply laced.}
Suppose that $C$ is nonsimply laced.
Then, in Step~2 above,
the factor  $M^{(a)}_1(u)$ in~\eqref{eq:ty4}
involves the term $T^{(b)}_{d_a}(u)$ for $a$ and $b$
with $a\sim b$, $d_a>1$, and $d_b=1$.
Therefore, Step~2 should be modif\/ied to def\/ine
these terms together.
For any~$Y$ satisfying~$\mathbb{Y}(C)$,
we construct a desired family $T$ in the following three steps:

{\em Step 1.} Choose arbitrarily
 $T^{(a)}_1(u)\in R^{\times}$
($a\in I$,
 $-\frac{d_a}{t}\leq  u
< \frac{d_a}{t}$).

{\em Step 2.}
Def\/ine
 $T^{(a)}_1(u)$
($a\in I$)
for the rest of $u\in\frac{1}{t}\mathbb{Z}$
as below.
(One can easily check that each step is well-def\/ined.)

Let $\{1<p_1<p_2<\cdots<p_k\}
=\{d_a\mid a\in I\}$,
and $I=I_1\sqcup I_{p_1}\sqcup I_{p_2}\sqcup \cdots
\sqcup I_{p_k}$, where $I_p:=\{a\in I\mid d_a=p\}$.

{\samepage

{\em Substep 1.}
\begin{itemize}\itemsep=0pt
\item[$(i)_1$.]
Def\/ine $T^{(a)}_1(-\frac{2}{t})$,
$T^{(a)}_1(\frac{1}{t})$ ($a\in I_1$) by
\begin{gather}
\label{eq:ty6}
{
T^{(a)}_{1}\left(u\pm \frac{d_a}{t}\right)
} =
\big(1+Y^{(a)}_1(u)^{-1}\big)
\frac{M^{(a)}_1(u)}
{T^{(a)}_{1}\left(u\mp \frac{d_a}{t}\right)}.
\end{gather}
\item[$(ii)_1$.]
Def\/ine $T^{(a)}_1(u)$
($a\in I_1$) for the rest of
$ -\frac{p_1}{t}\leq  u
< \frac{p_1}{t}$
by repeating $(i)_1$.
\end{itemize}}

{\em Substep 2.}
\begin{itemize}\itemsep=0pt
\item[$(i)_2$.]
Def\/ine $T^{(a)}_{p_1}(-\frac{1}{t})$,
$T^{(a)}_{p_1}(0)$ ($a\in I_1$) by
\begin{gather}
\label{eq:ty7}
T^{(a)}_{m+1}(u) =
\frac{1}{1+Y^{(a)}_m(u)}
\frac{T^{(a)}_m\left(u-\frac{d_a}{t}\right)
T^{(a)}_m\left(u+\frac{d_a}{t}\right)}
{T^{(a)}_{m-1}(u)},
\end{gather}
where $T^{(a)}_0(u)=1$.
\item[$(ii)_2$.]
Def\/ine $T^{(a)}_1(-\frac{p_1}{t}-\frac{1}{t})$,
$T^{(a)}_1(\frac{p_1}{t})$ ($a\in I_1\sqcup I_{p_1}$) by
\eqref{eq:ty6}.
\item[$(iii)_2$.]
Def\/ine $T^{(a)}_1(u)$
($a\in I_1\sqcup I_{p_1}$) for the rest of
$ -\frac{p_2}{t}\leq  u
< \frac{p_2}{t}$
by repeating $(i)_2$--$(ii)_2$.
\end{itemize}

{\em Substep 3.}
\begin{itemize}\itemsep=0pt
\item[$(i)_3$.]
Def\/ine $T^{(a)}_{p_1}(-\frac{p_2}{t}+\frac{p_1}{t}-\frac{1}{t})$,
$T^{(a)}_{p_1}(\frac{p_2}{t}-\frac{p_1}{t})$ ($a\in I_1$) by
\eqref{eq:ty7}.\\
Def\/ine $T^{(a)}_{p_2}(-\frac{1}{t})$,
$T^{(a)}_{p_2}(0)$ ($a\in I_1$) by \eqref{eq:ty7}.

\item[$(ii)_3$.]
Def\/ine $T^{(a)}_1(-\frac{p_2}{t}-\frac{1}{t})$,
$T^{(a)}_1(\frac{p_2}{t})$ ($a\in I_1\sqcup I_{p_1}\sqcup I_{p_2}$) by
\eqref{eq:ty6}.
\item[$(iii)_3$.]
Def\/ine $T^{(a)}_1(u)$
($a\in I_1\sqcup I_{p_1}\sqcup I_{p_2}$) for the rest of
$ -\frac{p_3}{t}\leq  u
< \frac{p_3}{t}$
by repeating $(i)_3$--$(ii)_3$.
\end{itemize}

{\em Substep 4.}
Continue to def\/ine $T^{(a)}_1(u)$ ($a\in I$) for the rest
of $u\in
\frac{1}{t} \mathbb{Z}$.

{\em Step 3.}
Def\/ine $T^{(a)}_m(u)$ ($a\in I, u\in \frac{1}{t}\mathbb{Z}$)
for $m\geq 2$
by \eqref{eq:ty7}.

Notice that
$T^{(a)}_1(u)$ is always def\/ined by \eqref{eq:ty6}, while
$T^{(a)}_m(u)$ for $m\geq 2$ is so by \eqref{eq:ty7}.
Then, the rest of the proof can be done in a  parallel way
to the simply laced case by using~\eqref{eq:Tidentity1} and~\eqref{eq:Tidentity2}.
\end{proof}

\begin{Definition}
\label{defn:Yring}
Let $\EuScript{Y}(C)$
be the commutative ring over $\mathbb{Z}$ with generators
$Y^{(a)}_m(u)^{\pm 1}$, \linebreak \mbox{$(1+Y^{(a)}_m(u))^{-1}$}
 ($a\in I$, $m\in \mathbb{N}$,
$u\in U $)
and the relations $\mathbb{Y}(C)$.
\end{Definition}

In terms of $\EuScript{T}(C)$ and $\EuScript{Y}(C)$,
 Theorem~\ref{thm:TtoY1} is rephrased as follows
(cf.~\cite[Theorem 2.12]{IIKNS} for the quantum af\/f\/ine algebras):

\begin{Theorem}\label{thm:TtoY2}\qquad

$(1)$ There is a ring homomorphism
\begin{gather*}
\varphi: \ \EuScript{Y}(C)  \rightarrow
\EuScript{T}(C)
\end{gather*}
defined by
\begin{gather*}
Y^{(a)}_m(u) \mapsto
\frac{M^{(a)}_m(u)}
{T^{(a)}_{m-1}(u)T^{(a)}_{m+1}(u)}.
\end{gather*}

$(2)$ There is a $($not unique$)$ ring homomorphism
\begin{gather*}
\psi: \
{\EuScript{T}}(C)
\rightarrow
\EuScript{Y}(C)
\end{gather*}
such that $\psi\circ\varphi = \mathrm{id}_{\EuScript{Y}(C)}$.
\end{Theorem}

There is another variation of Theorem \ref{thm:TtoY1}.
Let $\EuScript{T}^{\times}(C)$ (resp.\ $\EuScript{Y}^{\times}(C)$)
be the multiplicative subgroup
of all the invertible elements of
$\EuScript{T}(C)$ (resp.\ $\EuScript{Y}(C)$).
Clearly, $\EuScript{T}^{\times}(C)$
is generated by $T^{(a)}_m(u)$'s,
while $\EuScript{Y}^{\times}(C)$
is generated
by $Y^{(a)}_m(u)$'s and $1+Y^{(a)}_m(u)$'s.

\begin{Theorem}\label{thm:TtoY4}\qquad

$(1)$ There is a multiplicative group homomorphism
\begin{gather*}
\varphi: \  \EuScript{Y}^{\times}(C) \rightarrow\EuScript{T}^{\times}(C)
\end{gather*}
defined by
\begin{gather*}
Y^{(a)}_m(u) \mapsto
\frac{M^{(a)}_m(u)}
{T^{(a)}_{m-1}(u)T^{(a)}_{m+1}(u)},
\\
1+Y^{(a)}_m(u) \mapsto
\frac{T^{(a)}_{m}\bigl(u- \frac{d_a}{t}\bigr)
T^{(a)}_{m}\left(u+ \frac{d_a}{t}\right)}
{T^{(a)}_{m-1}(u)T^{(a)}_{m+1}(u)}.
\end{gather*}

$(2)$ There is a $($not unique$)$ multiplicative group homomorphism
\begin{gather*}
\psi: \
{\EuScript{T}}^{\times}(C)
\rightarrow
\EuScript{Y}^{\times}(C)
\end{gather*}
such that $\psi\circ\varphi = \mathrm{id}_{\EuScript{Y}(C)^{\times}}$.
\end{Theorem}

\section{Restricted T and Y-systems}
\label{sect:rest}

Here we introduce  a series of reductions of the
systems $\mathbb{T}(C)$ and $\mathbb{Y}(C)$
called the restricted T and Y-systems.
The restricted T and Y-systems for
the quantum af\/f\/ine algebras
are important
in application
to various integrable models.

We def\/ine integers $t_a$ ($a\in I$) by
\begin{gather*}
t_a=\frac{t}{d_a}.
\end{gather*}
\begin{Definition}
Fix an integer $\ell\geq 2$.
For a tamely laced generalized Cartan matrix $C$,
the {\em level $\ell$ restricted T-system $\mathbb{T}_{\ell}(C)$
 associated with $C$ (with the unit boundary condition)}
is the system of relations \eqref{eq:T1} and \eqref{eq:T2}
naturally restricted to
a family of variables
 $T_{\ell}=\{T^{(a)}_m(u) \mid a\in I; m=1,\dots,t_a\ell-1;
u\in U \}$,
where
$T^{(a)}_0 (u)= 1$, and furthermore, $T^{(a)}_{t_a\ell}(u)=1$
(the {\em unit boundary condition\/}) if they occur
in the right hand sides in the relations.
\end{Definition}

\begin{Definition}
Fix an integer $\ell\geq 2$.
For a tamely laced generalized Cartan matrix $C$,
the {\em level $\ell$ restricted Y-system $\mathbb{Y}_{\ell}(C)$
 associated with $C$}
is the system of relations \eqref{eq:Y1} and \eqref{eq:Y2}
naturally restricted to
a family of variables
 $Y_{\ell}=\{Y^{(a)}_m(u) \mid a\in I; m=1,\dots,t_a\ell-1;
u\in U \}$,
where
$Y^{(a)}_0 (u)^{-1}= 0$, and furthermore, $Y^{(a)}_{t_a\ell}(u)^{-1}=0$
if they occur
in the right hand sides in the relations.
\end{Definition}

The restricted version of Theorem
\ref{thm:TtoY1}~(1) holds.

\begin{Theorem}
\label{thm:TtoY3}
Let $R$ be any commutative ring over $\mathbb{Z}$ with identity element.
For any family
$T_{\ell}=\{T^{(a)}_m(u)\in R^{\times}
 \mid a\in I; m=1,\dots,t_a\ell-1; u\in U\}$
satisfying $\mathbb{T}_{\ell}(C)$ in $R$,
define a family
$Y_{\ell}=\{Y^{(a)}_m(u)\in R^{\times} \mid
 a\in I; m=1,\dots,t_a\ell-1; u\in U\}$ by
\eqref{eq:TtoY1}, where $T^{(a)}_0(u)=T^{(a)}_{t_a\ell}(u)=1$.
Then, $Y_{\ell}$ satisfies $\mathbb{Y}_{\ell}(C)$ in $R$.
\end{Theorem}
\begin{proof}

The calculation is formally the same as the one for
Theorem~\ref{thm:TtoY1}.
We have only to take care of the boundary term
\begin{gather}
\label{eq:bnd1}
\frac{1}{1+Y^{(a)}_{t_a\ell}(u)^{-1}}
=\frac{M^{(a)}_{t_a\ell}(u)}
{
T^{(a)}_{t_a\ell}
\bigl(u-\frac{d_a}{t}\bigr)
T^{(a)}_{t_a\ell}
\bigl(u+\frac{d_a}{t}\bigr)
},
\end{gather}
which formally appears in the right hand sides
of \eqref{eq:Y1} and \eqref{eq:Y2} for $m=t_a \ell-1$.
Since
\begin{gather*}
M^{(a)}_{t_a\ell}(u)=
\begin{cases}
\displaystyle
\prod_{b:b \sim a}
T^{(b)}_{t_b \ell}(u),
& d_a > 1,\\
\displaystyle
\prod_{b: b\sim a}
\left\{
\prod_{k=1}^{d_b}
T^{(b)}_{t_b\ell}
\left(u+\frac{1}{t}(2k-1-d_b)\right)
\right\},
& d_a = 1,
\end{cases}
\end{gather*}
the right hand side of~\eqref{eq:bnd1} is 1
under the boundary condition $T^{(a)}_{t_a\ell}(u)=1$ of
$\mathbb{T}_{\ell}(C)$.
This is compatible with $Y^{(a)}_{t_a\ell}(u)^{-1}=0$.
\end{proof}

Unfortunately the restricted version of Theorem
\ref{thm:TtoY1}~(2) does not hold due to
the boundary condition of~$\mathbb{T}_{\ell}(C)$.

\section{T and Y-systems from cluster algebras}
\label{sect:cluster}

In this section we introduce
T and Y-systems {\em associated with a class of cluster
algebras\/} \cite{FZ2,FZ4}
by generalizing some of the results
in \cite{FZ4,HL1,DK1,Kel,IIKNS,DK2}.
They include the restricted T and Y-systems of {\em simply laced type\/}
in Section \ref{sect:rest} as special cases.

\subsection[Systems $\mathbb{T}(B)$ and $\mathbb{Y}_{\pm}(B)$]{Systems $\boldsymbol{\mathbb{T}(B)}$ and $\boldsymbol{\mathbb{Y}_{\pm}(B)}$}

We warn the reader that the matrix $B$ in this section
 is dif\/ferent from the
one in Section~\ref{sect:qKM} and should not be confused.

\begin{Definition}[\cite{FZ1}]
An integer matrix
$B=(B_{ij})_{i,j\in I}$  is {\em skew-symmetrizable\/}
if there is a~diagonal matrix $D=\mathrm{diag}
(d_i)_{i\in I}$ with $d_i\in \mathbb{N}$
such that $DB$ is skew-symmetric.
For a skew-symmetrizable matrix $B$ and
$k\in I$, another matrix $B'=\mu_k(B)$,
called the {\em mutation of $B$ at $k$\/}, is def\/ined by
\begin{gather}
\label{eq:Bmut}
B'_{ij}=
\begin{cases}
-B_{ij}, & \mbox{$i=k$ or $j=k$},\\
B_{ij}+\frac{1}{2}
(|B_{ik}|B_{kj} + B_{ik}|B_{kj}|),
&\mbox{otherwise}.
\end{cases}
\end{gather}
\end{Definition}

The matrix $\mu_k(B)$ is also skew-symmetrizable.
The matrix mutation plays a central role in
the theory of cluster algebras.

We impose the following conditions on a skew-symmetrizable
matrix $B$:
The index set $I$ admits the decomposition
$I=I_+\sqcup I_-$ such that
\begin{gather}
\label{eq:B1}
\mbox{if $B_{ij}\neq 0$, then $(i,j)\in I_+\times I_-$
or $(i,j)\in I_-\times I_+$}.
\end{gather}
Furthermore, for composed mutations
$\mu_{+}=\prod_{i\in I_+} \mu_i$
and $\mu_{-}=\prod_{i\in I_-} \mu_i$,
\begin{gather}
\label{eq:B2}
\mu_+(B)=\mu_-(B)=-B.
\end{gather}

Note that $\mu_{\pm}(B)$ does not depend on the order of the product
due to \eqref{eq:B1}.

\begin{Lemma}
Under the condition \eqref{eq:B1}, the condition
\eqref{eq:B2} is equivalent to the following one:
For any $i,j\in I_+$,
\begin{gather}
\label{eq:BB}
\sum_{k: B_{ik}>0,B_{kj}>0}
B_{ik}B_{kj}
=
\sum_{k: B_{ik}<0,B_{kj}<0}
B_{ik}B_{kj}.
\end{gather}
The same holds for $i,j\in I_-$.
\end{Lemma}

\begin{proof}
Suppose that \eqref{eq:B2} holds.
Then, for any $i,j\in I_+$,
$\mu_-(B)_{ij}=-B_{ij}=0$
by~\eqref{eq:B1}.
It follows from~\eqref{eq:Bmut} that
\begin{gather*}
\sum_{k: B_{ik}B_{kj}>0} |B_{ik}|B_{kj}=0.
\end{gather*}
Therefore, we have \eqref{eq:BB}.
The rest of the proof is similar.
\end{proof}

\begin{Definition}
For a skew-symmetrizable matrix $B$ satisfying
the conditions \eqref{eq:B1} and \eqref{eq:B2},
the {\em T-system $\mathbb{T}(B)$
 associated with $B$}
is the following system of relations
for a family of variables
 $T=\{T_{i}(u) \mid i\in I,u\in\mathbb{Z} \}$:
\begin{gather*}
T_{i}(u-1)T_{i}(u+1)
=
\prod_{j: B_{ji}>0} T_j(u)^{B_{ji}}
+
\prod_{j: B_{ji}<0} T_j(u)^{-B_{ji}}.
\end{gather*}
\end{Definition}

For Y-systems, it is natural to introduce two kinds of systems.

\begin{Definition}
For  a skew-symmetrizable matrix $B$ satisfying
the conditions \eqref{eq:B1} and \eqref{eq:B2},
the {\em Y-systems $\mathbb{Y}_{+}(B)$ and $\mathbb{Y}_{-}(B)$
 associated with $B$}
are the following systems of relations
for a~family of variables
 $Y=\{Y_{i}(u) \mid i\in I,u\in \mathbb{Z} \}$,
respectively:
For $\mathbb{Y}_+(B)$,
\begin{gather}
\label{eq:Ybi+}
Y_{i}(u-1)Y_{i}(u+1)
=
\frac{
{
\prod\limits_{j: B_{ji}\gtrless0}
} (1+Y_j(u))^{\pm B_{ji}}}
{
{
\prod\limits_{j: B_{ji}\lessgtr0}
}
 (1+Y_j(u)^{-1})^{\mp B_{ji}}
},
\qquad
i\in I_{\pm}.
\end{gather}
For $\mathbb{Y}_-(B)$,
\begin{gather*}
Y_{i}(u-1)Y_{i}(u+1)
=
\frac{
{
\prod\limits_{j: B_{ji}\lessgtr0}
} (1+Y_j(u))^{\mp B_{ji}}}
{
{
\prod\limits_{j: B_{ji}\gtrless0}
}
 (1+Y_j(u)^{-1})^{\pm B_{ji}}
},
\qquad
i\in I_{\pm}.
\end{gather*}
\end{Definition}

\begin{Remark}
\label{rem:YY}
Two systems, $\mathbb{Y}_+(B)$ and $\mathbb{Y}_-(B)$, are transformed
into each other by any of the exchanges,
$Y_i(u)\leftrightarrow Y_i(u)^{-1}$,
$I_{\pm}\leftrightarrow I_{\mp}$,
or $B\leftrightarrow -B$.
There is no preferred choice
between~$\mathbb{Y}_+(B)$ and~$\mathbb{Y}_-(B)$, {\em a priori\/},
  and a convenient one can be used depending on the context.
On the contrary, $\mathbb{T}(B)$ is invariant under either
$I_{\pm}\leftrightarrow I_{\mp}$
or $B\leftrightarrow -B$.
\end{Remark}

\begin{Theorem}
\label{thm:TtoY5}
Let $R$ be any commutative ring over $\mathbb{Z}$ with identity element.
For any family
$T=\{T_i(u)\in R^{\times}
 \mid i\in I, u\in \mathbb{Z}\}$
satisfying $\mathbb{T}(B)$ in $R$,
define a family
$Y=\{Y_i(u)\in R^{\times} \mid
 i\in I,$ $u\in \mathbb{Z} \}$ by
\begin{gather*}
Y_i(u)=\prod_{j\in I}T_j(u)^{\pm B_{ji}},
\qquad
i\in I_{\pm}.
\end{gather*}
Then, $Y$ satisfies $\mathbb{Y}_{+}(B)$.
Similarly, define
 a family
$Y$ by
\begin{gather*}
Y_i(u)=\prod_{j\in I}T_j(u)^{\mp B_{ji}},
\qquad
i\in I_{\pm}.
\end{gather*}
Then, $Y$ satisfies $\mathbb{Y}_-(B)$ in $R$.
\end{Theorem}
\begin{proof}
By Remark \ref{rem:YY}, it is enough to prove the f\/irst statement only.
Then,
\begin{gather}
\label{eq:TtoY4}
Y_i(u) =\frac{
\prod\limits_{j:B_{ji}\gtrless 0}
T_j(u)^{\pm B_{ji}}
}
{
\prod\limits_{j:B_{ji}\lessgtr 0}
T_j(u)^{\mp B_{ji}}
},
\qquad i\in I_{\pm},
\\
\label{eq:TtoY5}
1+Y_i(u) =
\frac{
T_i(u-1)T_i(u+1)
}
{
\prod\limits_{j:B_{ji}\lessgtr 0}
T_j(u)^{\mp B_{ji}}
},
\qquad i\in I_{\pm},
\\
\label{eq:TtoY6}
1+Y_i(u)^{-1} =
\frac{
T_i(u-1)T_i(u+1)
}
{
\prod\limits_{j:B_{ji}\gtrless 0}
T_j(u)^{\pm B_{ji}}
},
\qquad i\in I_{\pm}.
\end{gather}
Note that for $j$ in the right hand side of \eqref{eq:Ybi+},
$j\in I_{\mp}$ by \eqref{eq:B1}.
By putting \eqref{eq:TtoY4}--\eqref{eq:TtoY6} into~\eqref{eq:Ybi+}, the right hand side of~\eqref{eq:Ybi+} is
\begin{gather*}
\prod_{j\in I}\bigl\{ T_j(u-1)T_j(u+1)
\bigr\}^{\pm B_{ji}}\!\!\prod_{j: B_{ji}\gtrless 0 }\!\Biggl\{
\prod_{k:B_{kj}\gtrless 0} T_k(u)^{\mp B_{kj}}
\Biggr\}^{\pm B_{ji}}\!\!\!
\prod_{j: B_{ji}\lessgtr 0 }\!\Biggl\{
\prod_{k:B_{kj}\lessgtr 0} T_k(u)^{\pm B_{kj}}
\Biggr\}^{\pm B_{ji}}\\
\qquad \overset{\eqref{eq:BB}}{=}
\prod_{j\in I}\bigl\{ T_j(u-1)T_j(u+1)
\bigr\}^{\pm B_{ji}},
\end{gather*}
which is the left hand side of \eqref{eq:Ybi+}.
\end{proof}

\subsection{Examples}
\label{subsec:Bexample}
Let us present some examples of
$\mathbb{T}(B)$ and $\mathbb{Y}_{\pm}(B)$.

\begin{Definition}
A symmetrizable generalized Cartan matrix
$C=(C_{ij})_{i,j\in I}$ is said to be
 {\em bipartite\/} if
the index set $I$ admits the decomposition
$I=I_+\sqcup I_-$ such that
\begin{gather*}
\mbox{if $C_{ij}<0$, then $(i,j)\in I_+\times I_-$
or $(i,j)\in I_-\times I_+$}.
\end{gather*}
\end{Definition}

\begin{Example}[\cite{FZ2,FZ4}]
\label{example:B1}
Let $C$ be a
{\em bipartite\/} symmetrizable generalized Cartan matrix,
which is not necessarily tamely laced.
Def\/ine the matrix $B=B(C)$ by
\begin{gather}
\label{eq:Bmat}
B_{ij}=
\begin{cases}
-C_{ij}, & (i,j)\in I_+ \times I_-,\\
C_{ij}, & (i,j)\in I_- \times I_+,\\
0, & \mbox{otherwise}.
\end{cases}
\end{gather}
The rule \eqref{eq:Bmat} is visualized in the diagram:
\begin{gather*}
\begin{matrix}
& _{ -C} &\\
+ &\rightarrow& -
\end{matrix}
\end{gather*}
Then, $B$ is skew-symmetrizable and
satisf\/ies the conditions \eqref{eq:B1} and \eqref{eq:B2}.
The correspon\-ding~$\mathbb{T}(B)$ and~$\mathbb{Y}_-(B)$ are given by
\begin{gather*}
T_i(u-1)T_i(u+1)
 =
1+ \prod_{j:j\sim i}
T_j(u)^{-C_{ji}},\\
Y_i(u-1)Y_i(u+1)
 =
\prod_{j:j\sim i}
(1+Y_j(u))^{-C_{ji}},
\end{gather*}
where $j\sim i$ means $C_{ji}<0$.
These systems are studied in \cite{FZ2,FZ4}.
When $C$ is {\em bipartite and simply laced},
they coincide with $\mathbb{T}_2(C)$ and $\mathbb{Y}_2(C)$
(for $U=\mathbb{Z}$) in Section \ref{sect:rest}.
When $C$ is {\em bipartite, tamely laced,
but nonsimply laced}, they are {\em different\/} from
$\mathbb{T}_2(C)$ and $\mathbb{Y}_2(C)$,
because the latter include factors depending on $u+\alpha$ ($\alpha\neq 0$)
in the right hand sides.
\end{Example}

\begin{Example}[{Square product
\cite{HL1,DK1,Kel,IIKNS,DK2}}]
\label{example:B2}
Let $C=(C_{ij})_{i,j\in I}$ and $C'=(C'_{i'j'})_{i',j'\in I'}$
 be a pair of {\em bipartite\/}
symmetrizable generalized Cartan matrices
with
$I=I_+\sqcup I_-$ and $I'=I'_+\sqcup I'_-$,
which are not necessarily tamely laced.
For $\mathbf{i}=(i,i')\in I\times I'$,
let us write $\mathbf{i}:(++)$ if $(i,i')\in I_+\times I'_+$, {\em etc}.
Def\/ine the matrix $B=(B_{\mathbf{i}\mathbf{j}})_{\mathbf{i},\mathbf{j}
\in I\times I'}$ by
\begin{gather}
\label{eq:Bsq}
B_{\mathbf{i}\mathbf{j}}=
\begin{cases}
-C_{ij}\delta_{i'j'},
 &
 \mathbf{i}:(-+), \ \mathbf{j}:(++)
\ \mbox{or}\
 \mathbf{i}:(+-), \ \mathbf{j}:(--),
\\
C_{ij}\delta_{i'j'},
 &
 \mathbf{i}:(++), \ \mathbf{j}:(-+)
\ \mbox{or} \
 \mathbf{i}:(--), \ \mathbf{j}:(+-),
\\
-\delta_{ij}C'_{i'j'},
 &
 \mathbf{i}:(++), \ \mathbf{j}:(+-)
\ \mbox{or} \
 \mathbf{i}:(--), \ \mathbf{j}:(-+),
\\
\delta_{ij}C'_{i'j'},
 &
 \mathbf{i}:(+-),  \ \mathbf{j}:(++)
\ \mbox{or} \
 \mathbf{i}:(-+), \ \mathbf{j}:(--),
\\
0, & \mbox{otherwise}.
\end{cases}
\end{gather}
The rule \eqref{eq:Bsq} is visualized in
the diagram:
\begin{gather}
\label{eq:square1}
\begin{matrix}
&& _{-C}&&\\
&(+-)& \rightarrow & (--)&\\
_{-C'}&\uparrow &&\downarrow& _{-C'}\\
&(++)&\leftarrow & (-+)&\\
&& _{-C}&&
\end{matrix}
\end{gather}
Since it generalizes the {\em square product of quivers\/} by \cite{Kel},
we call the matrix $B$ the {\em square product $B(C)\square B(C')$}
of the matrices $B(C)$ and $B(C')$ of \eqref{eq:Bmat}.

\begin{Lemma}
The matrix $B$ in \eqref{eq:Bsq} is skew-symmetrizable and
satisfies the conditions \eqref{eq:B1} and \eqref{eq:B2}
for $(I\times I')_+:=(I_+\times I'_+)\sqcup(I_-\times I'_-)$
and
$(I\times I')_-:=(I_+\times I'_-)\sqcup(I_-\times I'_+)$.
\end{Lemma}
\begin{proof}
Let $\mathrm{diag}(d_i)_{i\in I}$ and $\mathrm{diag}(d'_i)_{i\in I'}$ be
the diagonal  matrices  skew-symmetrizing $C$ and $C'$, respectively,
and let $D=\mathrm{diag}(d_id'_{i'})_{(i,i')\in
I\times I'}$.
Then, the matrix $DB$ is skew-symmetric.
The condition \eqref{eq:B1} is clear from \eqref{eq:square1}.
To show \eqref{eq:BB}, suppose, for example,
that $\mathbf{i}=(i,i'):(++)$ and
$\mathbf{j}=(j,j'):(--)$.
Then, $B_{\mathbf{i}\mathbf{k}}B_{\mathbf{k}\mathbf{j}}\neq 0$
only for $\mathbf{k}=(i,j')$ or $\mathbf{k}=(j,i')$;
furthermore,
$B_{\mathbf{i}\mathbf{k}},B_{\mathbf{k}\mathbf{j}}\geq 0$
(resp.\ $\leq 0$)
 for
$\mathbf{k}=(i,j')$ (resp.\ $\mathbf{k}=(j,i')$),
and $B_{\mathbf{i}\mathbf{k}}B_{\mathbf{k}\mathbf{j}}
=C_{ij}C'_{i'j'}$ for both.
Thus, \eqref{eq:BB} holds.
The other cases are similar.
\end{proof}

The corresponding $\mathbb{T}(B)$ and
$\mathbb{Y}_+(B)$ are given by
\begin{gather*}
T_{ii'}(u-1)T_{ii'}(u+1)
 =
\prod_{j:j\sim i}
T_{ji'}(u)^{-C_{ji}}
+
\prod_{j':j'\sim i'}
T_{ij'}(u)^{-C'_{j'i'}},\\
Y_{ii'}(u-1)Y_{ii'}(u+1)
 =
\frac
{
\prod\limits_{j:j\sim i}
(1+Y_{ji'}(u))^{-C_{ji}}
}
{
\prod\limits_{j':j'\sim i'}
(1+Y_{ij'}(u)^{-1})^{-C'_{j'i'}}
},
\end{gather*}
where $j\sim i$ and $j'\sim i'$ means $C_{ji}<0$
and $C'_{j'i'}<0$, respectively.
These systems slightly generalize the ones studied in
connection with
cluster algebras \cite{HL1,DK1,Kel,IIKNS,DK2}.
When  $C$ is bipartite and simply laced,
 and  $C'$ is the Cartan matrix of
type $A_{\ell-1}$
with $I'_+=\{1,3,\dots\}$ and $I'_-=\{2,4,\dots\}$,
$\mathbb{T}(B)$ and $\mathbb{Y}_+(B)$
 coincide with $\mathbb{T}_{\ell}(C)$ and $\mathbb{Y}_{\ell}(C)$
 in Section~\ref{sect:rest}. (The choice of~$I'_{\pm}$ is
not essential here.)
As in Example~\ref{example:B1},
when
$C$ is bipartite, tamely laced,
but nonsimply laced,
and $C'$ is the Cartan matrix of
type $A_{\ell-1}$,
they are {\em different\/} from
$\mathbb{T}_{\ell}(C)$ and $\mathbb{Y}_{\ell}(C)$.
\end{Example}

\begin{Example}
\label{example:B3}
Let us give an example which does not belong to
the classes in Examples \ref{example:B1} and \ref{example:B2}.
Let $B=(B_{ij})_{i,j\in I}$ with $I=\{1,\dots,7\}$ be the
skew-symmetric matrix whose positive components are given by
\begin{gather*}
B_{21} = B_{13}=2,\qquad
B_{34} =B_{35}=B_{36}=B_{37}=
B_{42}=B_{52}=B_{62}=B_{72}=1.
\end{gather*}
The matrix $B$ is represented by the following quiver:
\begin{gather*}
\includegraphics{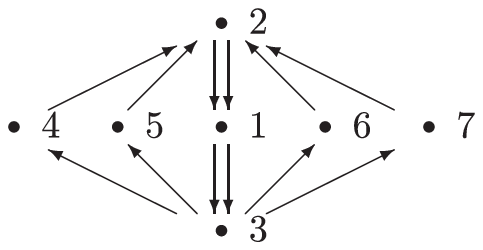}
\end{gather*}
With $I_+=\{2,3\}$ and $I_-=\{1,4,5,6,7\}$, the matrix $B$
satisf\/ies the conditions \eqref{eq:B1} and \eqref{eq:B2}.
\end{Example}

\subsection[$\mathbb{T}(B)$ and $\mathbb{Y}_{\pm}(B)$  as relations in cluster algebras]{$\boldsymbol{\mathbb{T}(B)}$ and $\boldsymbol{\mathbb{Y}_{\pm}(B)}$  as relations in cluster algebras}
\label{subsect:cluster}

The systems $\mathbb{T}(B)$ and $\mathbb{Y}_{\pm}(B)$ arise
as relations for {\em cluster variables\/}
and {\em coefficients\/}, respectively,
in the cluster algebra associated with $B$.
See \cite{FZ4,Kel} for def\/initions and information
for cluster algebras.

\subsubsection[$\mathbb{T}(B)$ and cluster algebras]{$\boldsymbol{\mathbb{T}(B)}$ and cluster algebras}

We start from T-systems.

\begin{Definition}
For a skew-symmetrizable matrix $B$ satisfying
the conditions \eqref{eq:B1} and \eqref{eq:B2},
let~$\EuScript{T}(B)$
be the commutative ring over $\mathbb{Z}$  with generators
$T_i(u)^{\pm 1}$ ($i\in I, u\in \mathbb{Z} $)
and the relations~$\mathbb{T}(B)$.
Also, let $\EuScript{T}^{\circ}(B)$
be the subring of $\EuScript{T}(B)$
generated by
$T_i(u)$ ($i\in I,u\in \mathbb{Z} $).
\end{Definition}

Let $\varepsilon:I \rightarrow \{+,-\}$ be the sign
function def\/ined by $\varepsilon(i)=\varepsilon$ for $i\in I_{\varepsilon}$.
For $(i,u)\in I\times \mathbb{Z}$,
we set the `parity conditions' $\mathbf{P}_{+}$ and
$\mathbf{P}_{-}$ by
\begin{gather*}
\mathbf{P}_{\pm}: \
\varepsilon(i)(-1)^{u} = \pm,
\end{gather*}
where we identify $+$ and $-$ with $1$ and $-1$, respectively.
For $\varepsilon\in \{+,-\}$,
def\/ine $\EuScript{T}^{\circ}(B)_{\varepsilon}$
to be the subring of $\EuScript{T}^{\circ}(B)$
generated by
those $T_i(u)$ with $(i,u)$ satisfying
$\mathbf{P}_{\varepsilon}$.
Then, we have $
\EuScript{T}^{\circ}(B)_+
\simeq
\EuScript{T}^{\circ}(B)_-
$ by $T_i(u)\mapsto T_i(u+1)$ and
\begin{gather*}
\EuScript{T}^{\circ}(B)
\simeq
\EuScript{T}^{\circ}(B)_+
\otimes_{\mathbb{Z}}
\EuScript{T}^{\circ}(B)_-.
\end{gather*}

Let
$\mathcal{A}(B,x)$
be the {\em cluster algebra
 with trivial coefficients},
where $(B,x)$ is the initial seed \cite{FZ4}.
We set $x(0)=x$ and def\/ine clusters
 $x(u)=(x_i(u))_{i\in I}$ ($u\in \mathbb{Z}$)
by the sequence of mutations
\begin{gather}
\label{eq:QseqADE}
\cdots
 \overset{\mu_-}{\longleftrightarrow}
(B,x(0))
\overset{\mu_+}{\longleftrightarrow}
(-B,x(1))
\overset{\mu_-}{\longleftrightarrow}
(B,x(2))
\overset{\mu_+}{\longleftrightarrow}
\cdots.
\end{gather}

\begin{Definition}The {\em T-subalgebra
${\mathcal{A}}_T(B,x)$
 of
${\mathcal{A}}(B,x)$
associated with
the sequence \eqref{eq:QseqADE}}
is the subring of
$\mathcal{A}(B,x)$
generated by $x_i(u)$ ($i\in I$, $u\in \mathbb{Z}$).
\end{Definition}
The ring ${\mathcal{A}}_T(B,x)$
 is no longer a cluster algebra in general,
because it is not closed under mutations.

\begin{Lemma}\label{lem:x} \quad

$(1)$ $x_i(u)=x_i(u\mp1)$ for $(i,u)$ satisfying $\mathbf{P}_{\pm}$.

$(2)$ The family $x=\{ x_i(u)\mid
i\in I, u\in \mathbb{Z}
\}$ satisfies the T-system $\mathbb{T}(B)$
 by replacing~$T_i(u)$
 in~$\mathbb{T}(B)$ with~$x_i(u)$.
\end{Lemma}

\begin{proof}
This follows from the exchange relation
of a cluster $x$ by the
 mutation $\mu_k$ \cite{FZ4}:
\begin{gather*}
x'_i =
\begin{cases}
{x_i},&i\neq k,\\
\displaystyle \frac{1}{x_k}
\left(\prod_{j:B_{ji}>0} x^{B{_{ji}}}
+
\prod_{j:B_{ji}<0} x^{-B{_{ji}}}\right),
& i = k.
\end{cases}\tag*{\qed}
\end{gather*}
\renewcommand{\qed}{}
\end{proof}

\begin{Theorem}[cf.\ \text{\cite[Proposition 4.24]{IIKNS}}]
\label{thm:TAA}
For $\varepsilon\in \{+,-\}$,
the ring $\EuScript{T}^{\circ}(B)_{\varepsilon}$
is isomorphic to
$
{\mathcal{A}}_T(B,x)
$
by the correspondence $T_i(u)\rightarrow x_i(u)$
for $(i,u)$ satisfying $\mathbf{P}_{\varepsilon}$.
\end{Theorem}
\begin{proof}
This follows from Lemma \ref{lem:x}
using the same argument  as the one for \cite[Proposition~4.2]{IIKNS}.
\end{proof}

\subsubsection[$\mathbb{Y}_{\varepsilon}(B)$ and cluster algebras]{$\boldsymbol{\mathbb{Y}_{\varepsilon}(B)}$ and cluster algebras}

We present a parallel result for Y-systems.

A {\em semifield\/} $(\mathbb{P},+)$ is an
abelian multiplicative group $\mathbb{P}$ endowed with a binary
operation of addition $+$ which is commutative,
associative, and distributive with respect to the
multiplication in $\mathbb{P}$ \cite{FZ4,HW}.
(Here we use the symbol $+$ instead of $\oplus$ in \cite{FZ4}
to make the description a little simpler.)

\begin{Definition}
\label{def:YB}
For $\varepsilon\in \{+,-\}$
and a skew-symmetrizable matrix $B$ satisfying
the conditions~\eqref{eq:B1} and~\eqref{eq:B2},
 let $\tilde{\EuScript{Y}}_{\varepsilon}(B)$
be the semif\/ield with generators
$Y_i(u)$  ($i\in I$, $u\in \mathbb{Z} $)
and the relations~$\mathbb{Y}_{\varepsilon}(B)$.
Let $\tilde{\EuScript{Y}}^{\circ}_{\varepsilon}(B)$
be the multiplicative subgroup
of $\tilde{\EuScript{Y}}_{\varepsilon}(B)$
generated by
$Y_i(u)$  and $1+Y_i(u)$
 ($i\in I$, $u\in \mathbb{Z} $).
(We use the notation $\tilde{\EuScript{Y}}$ to distinguish
it from the ring $\EuScript{Y}$ in Def\/inition~\ref{defn:Yring}.)
\end{Definition}

Def\/ine $\tilde{\EuScript{Y}}^{\circ}_{\varepsilon}(B)_{+}$
(resp.\ $\tilde{\EuScript{Y}}^{\circ}_{\varepsilon}(B)_{-}$)
to be the subgroup of $\tilde{\EuScript{Y}}^{\circ}_{\varepsilon}(B)$
generated by
those $Y_i(u)$ and $1+Y_i(u)$
with  $(i,u)$ satisfying
$\mathbf{P}_+$ (resp.\ $\mathbf{P}_-)$.
Then, we have
$\tilde{\EuScript{Y}}^{\circ}_{\varepsilon}(B)_+
\simeq
\tilde{\EuScript{Y}}^{\circ}_{\varepsilon}(B)_-
$
by $Y_i(u)\mapsto Y_i(u+1)$ and
\begin{gather*}
\tilde{\EuScript{Y}}^{\circ}_{\varepsilon}(B)
\simeq
\tilde{\EuScript{Y}}^{\circ}_{\varepsilon}(B)_+
\times
\tilde{\EuScript{Y}}^{\circ}_{\varepsilon}(B)_-.
\end{gather*}

Let $\mathcal{A}(B,x,y)$
be the {\em cluster algebra
 with coefficients
in the universal
semifield
$\mathbb{Q}_{\mathrm{sf}}(y)$},
where $(B,x,y)$ is the initial seed  \cite{FZ4}.
To make the setting parallel to T-systems,
we introduce the {\em coefficient group $\mathcal{G}(B,y)$
associated with $\mathcal{A}(B,x,y)$},
 which is the multiplicative subgroup of
the semif\/ield $\mathbb{Q}_{\mathrm{sf}}(y)$ generated by all
the elements $y_i'$ of coef\/f\/icient tuples of $\mathcal{A}(B,x,y)$
together with $1+y_i'$.

We set $x(0)=x$, $y(0)=y$ and def\/ine
clusters $x(u)=(x_i(u))_{i\in I}$
 and coef\/f\/icient tuples $y(u)=(y_i(u))_{i\in I}$
 ($u\in \mathbb{Z}$)
by the sequence of mutations
\begin{gather}
\label{eq:QseqADE2}
\cdots
 \overset{\mu_-}{\longleftrightarrow}
(B,x(0),y(0))
\overset{\mu_+}{\longleftrightarrow}
(-B,x(1),y(1))
 \overset{\mu_-}{\longleftrightarrow}
(B,x(2),y(2))
\overset{\mu_+}{\longleftrightarrow}
\cdots.
\end{gather}

\begin{Definition}
The {\em Y-subgroup
${\mathcal{G}}_Y(B,y)$
of ${\mathcal{G}}(B,y)$
associated with
the sequence~\eqref{eq:QseqADE2}}
is the multiplicative subgroup of
${\mathcal{G}}(B,y)$
generated by
$y_i(u)$  and  $1+ y_i(u)$
($i\in I$, $u\in \mathbb{Z} $).
\end{Definition}

\begin{Lemma}\label{lem:y}\qquad

$(1)$ $y_i(u)=y_i(u\pm1)^{-1}$ for $(i,u)$ satisfying $\mathbf{P}_{\pm}$.

$(2)$ For $\varepsilon\in \{+,-\}$,
the family $y_{\varepsilon}=\{ y_i(u)\mid
\mbox{$(i,u)$ satisfying $\mathbf{P}_{\varepsilon}$}
\}$ satisf\/ies the Y-system $\mathbb{Y}_{\varepsilon}(B)$
 by replacing $Y_i(u)$
in $\mathbb{Y}_{\varepsilon}(B)$ with $y_i(u)$.
\end{Lemma}

\begin{proof}
This follows from the exchange relation
of a  coef\/f\/icient tuple $y$ by the
 mutation $\mu_k$ \cite{FZ4}:
\begin{gather*}
y'_i =
\begin{cases}
{y_k}{}^{-1}, &i=k,\\
y_i (1+ {y_k}{}^{-1})^{-B_{ki}}, &
i\neq k,\ B_{ki}\geq 0,\\
y_i (1+ y_k)^{-B_{ki}}, &
i\neq k,\ B_{ki}\leq 0.
\end{cases}\tag*{\qed}
\end{gather*}
\renewcommand{\qed}{}
\end{proof}

\begin{Theorem}
\label{thm:YAA}
The group $\tilde{\EuScript{Y}}^{\circ}_+(B)_{\pm}$
is isomorphic to
${\mathcal{G}}_Y(B,y)
$
by the correspondence $Y_i(u)\mapsto y_i(u)^{\pm1}$,
$1+Y_i(u)\mapsto 1+ y_i(u)^{\pm1}$
for $(i,u)$ satisfying $\mathbf{P}_{\pm}$.
Similarly, the group $\tilde{\EuScript{Y}}^{\circ}_-(B)_{\pm}$
is isomorphic to
${\mathcal{G}}_Y(B,y)$
by the correspondence $Y_i(u)\mapsto y_i(u)^{\mp1}$,
$1+Y_i(u)\mapsto 1+ y_i(u)^{\mp1}$ for
$(i,u)$ satisfying $\mathbf{P}_{\pm}$.
\end{Theorem}
\begin{proof}
Let us show that $\tilde{\EuScript{Y}}^{\circ}_{+}(B)_+
\simeq {\mathcal{G}}_Y(B,y)$.
Let
$f: \mathbb{Q}_{\mathrm{sf}}(y)
\rightarrow \tilde{\EuScript{Y}}_+(B)
$
be  the semif\/ield homomorphism def\/ined by
\begin{gather*}
f: \ y_i   \mapsto
\begin{cases}
Y_i(0), & i\in I_+,\\
Y_i(-1)^{-1}, & i\in I_-.\\
\end{cases}
\end{gather*}
Then,
 due to Lemma  \ref{lem:y}~(2),
it can be shown by induction on $\pm u$ that
we have
$f:y_i(u)\mapsto
Y_i(u)$
for any $(i,u)$ satisfying $\mathbf{P}_+$,
and
$f:y_i(u)\mapsto
Y_i(u-1)^{-1}$ for any $(i,u)$ satisfying $\mathbf{P}_-$.
By the restriction of $f$,
we have
a multiplicative group homomorphism
 $f': {\mathcal{G}}_Y(B,y)
\rightarrow \tilde{\EuScript{Y}}_+^{\circ}(B)_+$.
On the other hand, by Lemma  \ref{lem:y} (2) again,
a semif\/ield homomorphism
 $g:
\tilde{\EuScript{Y}}_{+}(B)
\rightarrow
\mathbb{Q}_{\mathrm{sf}}(y)
$ is def\/ined by
$ Y_i(u)
\mapsto y_i(u)^{\pm1}$ for $(i,u)$ satisfying $\mathbf{P}_{\pm}$.
By the restriction of $g$,
we have a multiplicative group homomorphism
 $g':
\tilde{\EuScript{Y}}^{\circ}_{+}(B)_+
\rightarrow
{\mathcal{G}}_Y(B,y)$.
Then, $f'$ and $g'$ are the inverse to each other
by Lemma \ref{lem:y} (1).
Therefore, $\tilde{\EuScript{Y}}^{\circ}_{+}(B)_+
\simeq
{\mathcal{G}}_Y(B,y)$.
The other cases are similar.
\end{proof}

\subsection{Restricted T and Y-systems and cluster algebras: simply laced case}

The
restricted T and Y-systems, $T_{\ell}(C)$ and $Y_{\ell}(C)$,
introduced in Section \ref{sect:rest}
are special cases of~$T(B)$ and $Y_{\pm}(B)$, if $C$ is simply laced.
Therefore, they are also related to cluster algebras.

\subsubsection{Bipartite case}
Suppose that $C$ is a {\em simply laced and bipartite\/}
 generalized Cartan matrix.
Then, we have already seen
in Examples \ref{example:B1} and
\ref{example:B2}
 that
$T_{\ell}(C)$ and $Y_{\ell}(C)$ coincides with
$T(B)$ and $Y_{\varepsilon}(B)$ for some $B$ and $\varepsilon$.
Therefore, we immediately obtain the following results as special cases of
Theorems \ref{thm:TAA} and \ref{thm:YAA}.

\begin{Corollary}
\label{cor:cluster}
Let $C$ be a  simply laced and bipartite
 generalized Cartan matrix
with $I=I_+\sqcup I_-$.
For $\varepsilon\in \{+,-\}$,
let $\EuScript{T}^{\circ}_{\ell}(C)_{\varepsilon}$
be the  subring of $\EuScript{T}^{\circ}_{\ell}(C)$ for $U=\mathbb{Z}$
 generated by $T^{(a)}_m(u)$ ($a\in I;m=1,\dots,\ell-1;u\in \mathbb{Z}$)
 satisfying $\varepsilon(a)(-1)^{m+1+u}
 = \varepsilon$.
Then, we have the following:

$(1)$ $\EuScript{T}^{\circ}_2(C)_{\varepsilon}$ is isomorphic to
${\mathcal{A}}_T(B,x)$ with $B=B(C)$
by the correspondence
$T^{(a)}_1(u)\mapsto x_a(u)$.

$(2)$ For $\ell \geq 3$,
$\EuScript{T}^{\circ}_{\ell}(C)_{\varepsilon}$ is isomorphic to
${\mathcal{A}}_T(B,x)$ with $B=B(C)\square B(C')$
by the correspondence
$T^{(a)}_m(u)\mapsto x_{am}(u)$, where
$C'$ is the Cartan matrix of type $A_{\ell-1}$ with
$I'_+=\{1,3,\dots\}$ and $I'_-=\{2,4,\dots\}$.
\end{Corollary}

\begin{Corollary}
\label{cor:cluster2}
Let $C$ be a  simply laced and bipartite
 generalized Cartan matrix
with $I=I_+\sqcup I_-$.
Let $\tilde{\EuScript{Y}}_{\ell}(C)$
be the  semifield with generators
$Y^{(a)}_m(u)$ $(a\in I;m=1,\dots,\ell-1;u\in \mathbb{Z})$
and the relations $\mathbb{Y}_{\ell}(C)$ for $U=\mathbb{Z}$.
Let $\tilde{\EuScript{Y}}^{\circ}_{\ell}(C)$ be the multiplicative
subgroup of $\tilde{\EuScript{Y}}_{\ell}(C)$
generated by~$Y^{(a)}_m(u)$ and $1+Y^{(a)}_m(u)$
 $(a\in I;m=1,\dots,\ell-1;u\in \mathbb{Z})$.
For $\varepsilon\in \{+,-\}$,
let $\EuScript{Y}^{\circ}_{\ell}(C)_{\varepsilon}$
be the  multiplicative
subgroup of $\tilde{\EuScript{Y}}^{\circ}_{\ell}(C)$
 generated by $Y^{(a)}_m(u)$ and $1+Y^{(a)}_m(u)$
 $(a\in I;m=1,\dots,\ell-1;$ $u\in \mathbb{Z})$
 satisfying $\varepsilon(a)(-1)^{m+1+u}
 = \varepsilon$.
Then, we have the following:

$(1)$ $\tilde{\EuScript{Y}}^{\circ}_2(C)_{\pm}$ is isomorphic to
${\mathcal{G}}_Y(B,y)$ with $B=B(C)$
by the correspondence
$Y^{(a)}_1(u)\mapsto y_a(u)^{\mp1}$,
$1+Y^{(a)}_1(u)\mapsto 1+y_a(u)^{\mp1}$.

$(2)$ For $\ell \geq 3$,
$\tilde{\EuScript{Y}}^{\circ}_{\ell}(C)_{\pm}$ is isomorphic to
${\mathcal{G}}_Y(B,y)$ with $B=B(C)\square B(C')$
by the correspondence
$Y^{(a)}_m(u)\mapsto y_{am}(u)^{\pm1}$,
$1+Y^{(a)}_m(u)\mapsto 1+y_{am}(u)^{\pm1}$,
where
$C'$ is the Cartan matrix of type $A_{\ell-1}$ with
$I'_+=\{1,3,\dots\}$ and $I'_-=\{2,4,\dots\}$.
\end{Corollary}

The slight discrepancy of the signs between $\ell=2$ and
$\ell\geq 3$ is due to the convention adopted here and
not an essential problem.

\subsubsection{Nonbipartite case}

Let us extend Corollaries \ref{cor:cluster}
and  \ref{cor:cluster2}
to a {\em simply laced and nonbipartite\/}
generalized Cartan matrix $C$.
The Cartan matrix of type $A^{(1)}_{2r}$ is
such an example.
In general, a generalized Cartan matrix $C$ is
bipartite if and only if there is no odd cycle
in the corresponding Dynkin diagram.
Without loss of generality we can assume that
$C$ is {\em indecomposable}; namely, the corresponding
Dynkin diagram is connected.

\begin{Definition}
Let $C=(C_{ij})_{i,j\in I}$  be a simply laced, nonbipartite, and
indecomposable generalized Cartan matrix.
We introduce an index set
$I^{\#}=I^{\#}_+\sqcup I^{\#}_-$,
where $I^{\#}_+=\{i_+\}_{i\in I}$ and
$I^{\#}_-=\{i_-\}_{i\in I}$,
and def\/ine a matrix $C^{\#}=(C^{\#}_{\alpha\beta})_{\alpha,\beta
\in I^{\#}}$ by
\begin{gather*}
C^{\#}_{\alpha\beta}=
\begin{cases}
2, & \alpha=\beta,\\
C_{ij},&\mbox{$(\alpha,\beta)=(i_+, j_-)$
or $(i_-,j_+)$},\\
0, & \mbox{otherwise.}
\end{cases}
\end{gather*}
We call $C^{\#}$  the {\em bipartite double\/} of $C$.
\end{Definition}

It is clear that $C^{\#}$ is a simply laced
and indecomposable generalized Cartan matrix;
furthermore, it is bipartite with
$I^{\#}=I^{\#}_+\sqcup I^{\#}_-$.

\begin{Example}
Let $C$ be the Cartan matrix corresponding to the Dynkin diagram
in the left hand side below.
Then,
$C^{\#}$ is the  Cartan matrix corresponding to the Dynkin diagram
in the right hand side.\vspace*{-2mm}
\begin{gather*}
\includegraphics{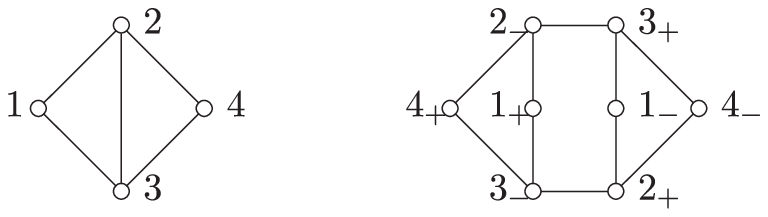}
\end{gather*}

\vspace*{-3mm}

\noindent
Here is another example.\vspace*{-2mm}
\begin{gather*}
\includegraphics{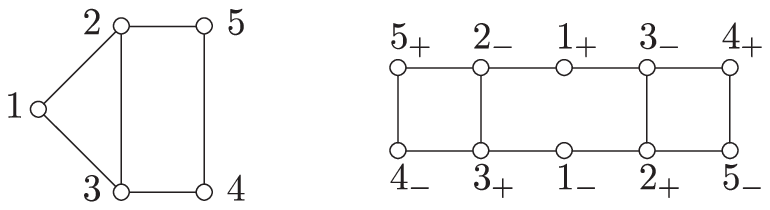}
\end{gather*}
\end{Example}

\begin{Proposition}
\label{prop:double}
Let $C=(C_{ij})_{i,j\in I}$ be a simply laced,
 nonbipartite, and indecomposable
generalized Cartan matrix, and $C^{\#}$ be its
bipartite double.

$(1)$
Let $\EuScript{T}^{\circ}_{\ell}(C^{\#})_+$ be the ring defined
in Corollary $\ref{cor:cluster}$.
Then, $\EuScript{T}^{\circ}_{\ell}(C)$ is isomorphic to
$\EuScript{T}^{\circ}_{\ell}(C^{\#})_+$
 by the correspondence
$T^{(a)}_m(u)
\mapsto
T^{(a_{\pm})}_m(u)$
for $(-1)^{m+1+u}=\pm$.

$(2)$
Let $\tilde{\EuScript{Y}}^{\circ}_{\ell}(C^{\#})_+$ be the
multiplicative group defined
in Corollary $\ref{cor:cluster2}$.
Then, $\tilde{\EuScript{Y}}^{\circ}_{\ell}(C)$ is isomorphic to
$\tilde{\EuScript{Y}}^{\circ}_{\ell}(C^{\#})_+$ by the correspondence
$Y^{(a)}_m(u)
\mapsto
Y^{(a_{\pm})}_m(u)$,
$1+Y^{(a)}_m(u)
\mapsto
1+Y^{(a_{\pm})}_m(u)$
for $(-1)^{m+1+u}=\pm$.
\end{Proposition}

\begin{proof}
The generators and relations
of the both sides coincide
under the correspondence.
\end{proof}

Combining Corollaries \ref{cor:cluster},
\ref{cor:cluster2}, and
Proposition \ref{prop:double},
we have the versions of Corolla\-ries~\ref{cor:cluster}
and~\ref{cor:cluster2}
in the nonbipartite case.

\begin{Corollary}
\label{cor:cluster3}
Let $C$ and $C^{\#}$ be the same ones as in
Proposition $\ref{prop:double}$.

$(1)$
$\EuScript{T}^{\circ}_2(C)$ is isomorphic to
${\mathcal{A}}_T(B,x)$ with $B=B(C^{\#})$
by the correspondence
$T^{(a)}_1(u)\mapsto x_{a_{\pm}}(u)$
for $(-1)^u = \pm$.

$(2)$ For $\ell \geq 3$,
$\EuScript{T}^{\circ}_{\ell}(C)$ is isomorphic to
${\mathcal{A}}_T(B,x)$ with $B=B(C^{\#})\square B(C')$
by the correspondence
$T^{(a)}_m(u)\mapsto x_{a_{\pm},m}(u)$
for $(-1)^{m+1+u} = \pm$,
where
$C'$ is the Cartan matrix of type $A_{\ell-1}$ with
$I'_+=\{1,3,\dots\}$ and $I'_-=\{2,4,\dots\}$.
\end{Corollary}

\begin{Corollary}
\label{cor:cluster4}
Let $C$ and $C^{\#}$ be the same ones as in
Proposition $\ref{prop:double}$.

$(1)$
$\tilde{\EuScript{Y}}^{\circ}_2(C)$ is isomorphic to
${\mathcal{G}}_Y(B,y)$ with $B=B(C^{\#})$
by the correspondence
$Y^{(a)}_1(u)\mapsto y_{a_{\pm}}(u)^{-1}$,
$1+Y^{(a)}_1(u)\mapsto 1+y_{a_{\pm}}(u)^{-1}$
for $(-1)^u = \pm$.

$(2)$ For $\ell \geq 3$,
$\tilde{\EuScript{Y}}^{\circ}_{\ell}(C)$ is isomorphic to
${\mathcal{G}}_Y(B,y)$ with $B=B(C^{\#})\square B(C')$
by the correspondence
$Y^{(a)}_m(u)\mapsto y_{a_{\pm},m}(u)$,
$1+Y^{(a)}_m(u)\mapsto 1+y_{a_{\pm},m}(u)$
for $(-1)^{m+1+u} = \pm$,
where
$C'$ is the Cartan matrix of type $A_{\ell-1}$ with
$I'_+=\{1,3,\dots\}$ and $I'_-=\{2,4,\dots\}$.
\end{Corollary}

\subsection{Concluding remarks}

One can further extend
Corollaries \ref{cor:cluster},
\ref{cor:cluster2},
\ref{cor:cluster3},
and \ref{cor:cluster4}
to the {\em tamely laced and nonsimply laced} case
by introducing T and Y-systems associated with
another class of cluster algebras\footnote{Inoue R., Iyama O., Keller B., Kuniba A., Nakanishi T., In preparation.}. 
Therefore, we conclude that
all the restricted T and Y-systems
associated with tamely laced generalized Cartan matrices
introduced in Section \ref{sect:rest}
are identif\/ied with the T and Y-systems
 associated with a certain class of cluster algebras.

The following question is left as an important problem:
What are the T and Y-systems associated with
{\em nontamely laced} symmetrizable generalized Cartan matrices?

\subsection*{Acknowledgements}
It is our great pleasure to thank
Professor Tetsuji Miwa
on the occasion of his sixtieth birthday
 for his generous
support
 and continuous interest in our works
through many years.
We also thank David Hernandez for a valuable comment.

\pdfbookmark[1]{References}{ref}
\LastPageEnding

\end{document}